\documentclass[preprint,3p,12pt]{elsarticle}

\usepackage{moreverb}
\usepackage{graphicx}
\usepackage{epstopdf}
\usepackage{epic}
\usepackage{ulem}

\usepackage{amssymb,amsfonts,amsmath,amscd}
\usepackage{makecell,multirow}
\usepackage{booktabs}

\usepackage{tikz}
\usetikzlibrary{arrows,shapes,chains}
\usetikzlibrary{positioning}
\usetikzlibrary{calc,3d}

\usepackage{algorithm}
\usepackage{mathrsfs}
\usepackage{algpseudocode}

\usepackage{theorem}
\newtheorem{remark}{Remark}

\usepackage{bigstrut}
\usepackage{subfig}
\usepackage{stfloats}

\newcommand{\red}[1]{{#1}}

\usepackage{lineno,hyperref}

\journal{J. Comput. Phys.}

\begin{document}
\captionsetup[figure]{labelfont={bf},name={Fig.},labelsep=period}
\begin{frontmatter}



\title{Parallel multilevel restricted Schwarz preconditioners for implicit simulation of subsurface flows with Peng-Robinson equation of state}
\author[label1]{Rui Li}
\ead{lirui319@hnu.edu.cn}
\author[label1]{Haijian Yang\corref{cor1}}
\ead{haijianyang@hnu.edu.cn}
\address[label1]{School of Mathematics, Hunan University, Changsha, Hunan 410082, PR China}
\author[label2]{Chao Yang}
\ead{chao\_yang@pku.edu.cn}
\address[label2]{School of Mathematical Sciences, Peking University, Beijing 100871, PR China}
\cortext[cor1]{Corresponding author. }

\begin{abstract}
Parallel algorithms and simulators with good scalabilities are particularly important for  large-scale reservoir simulations on modern  supercomputers with a large number of processors. In this paper, we introduce and study a family of highly scalable multilevel restricted additive Schwarz (RAS) methods for the fully implicit solution of subsurface flows with Peng-Robinson equation of state in two and three dimensions.  With the use of a second-order fully implicit scheme, the proposed simulator is unconditionally stable with the  relaxation of the time step size by the stability condition. The investigation then focuses on the development of several types of multilevel overlapping additive Schwarz methods for the preconditioning of the resultant linear system arising from the inexact Newton iteration,  and some fast solver technologies are presented for the assurance of the multilevel approach efficiency and scalability.  We numerically show that the proposed fully implicit framework is highly efficient for solving both standard benchmarks as well as realistic problems  with several hundreds of millions of unknowns and scalable to 8192 processors on the Tianhe-2 supercomputer.
\end{abstract}


\end{frontmatter}

\section{Introduction}\label{sc-int}
Simulation of subsurface fluid flows in porous media is currently an important topic of interest in many applications, such as hydrology and groundwater flow, oil and gas reservoir simulation, CO$_2$ sequestration,  and waste management  \cite{Chen,Firoozabadi}. The extensive growing demand on accurate modeling subsurface systems has produced persistent requirements on the understanding of many different physical processes, including multiphase fluid flow and transport,  fluid phase behavior with sophisticated equation of state, and geomechanical deformations.  To simulate these processes, one needs to solve a set of coupled, nonlinear, time-dependent partial differential equations (PDEs) with complicated nonlinear behaviors arising from the complexity of  the geological media and reservoir properties.  Additional computational challenge comes from  the involving of Peng-Robinson equation of state,  which governs the equilibrium distribution of fluid and has a remarkable influence on the high variation of the Darcy velocity. Due to their high computational complexity, numerically approximating subsurface phenomena with high resolution is critical to the reservoir engineering  for accurate predictions of costly projects. Hence, parallel reservoir simulators equipped with  robust and scalable solvers  on high performance computing platforms  are crucial to achieve efficient simulations of this intricate problem \cite{Chen07}.

In terms of the numerical solution of partial differential equations, the fully implicit  method  \cite{bui17,bui18,Lacroix,Liu2018,qiao,Sko13,wang,haijian-cmame} is a widely preferred approach to solve various problems arising {from the discretization of subsurface fluid flows}.  The fully implicit scheme is unconditionally stable and can relax the stability requirement of the Courant-Friedrichs-Lewy (CFL) condition. It therefore provides a consistent and robust way to simulate the subsurface fluid flow problem  in long-term and extreme-scale computations.  In particular, a great advantage of the fully implicit algorithm is that the corresponding nonlinear equations are  implicitly solved  in a monolithic way. This characteristic feature strengthens its potential to allow the  addition of more physics and the introduction of more equations without changing much of the simulator framework, which greatly expands the scope of application of the fully implicit approach. In spite of being stable for arbitrary large time steps, when a fully implicit scheme is applied, one must solve a large, sparse, and ill-conditioned linear system at each nonlinear iteration. It remains challenging and important to design robust and scalable solvers for the large scale simulations on high performance computing platforms.  In this study, we employ the framework of Newton--Krylov algorithms \cite{Cai98,cs02a,fim5,Knoll04,liu,haijian-sisc1} to guarantee the nonlinear consistency, and mainly focus our efforts on designing an efficient preconditioning strategy to substantially reduce the condition number of the linear system.

In reservoir simulations, some effective preconditioning techniques, such as block preconditioners \cite{Haga,Lee,White}, the Constrained Residual Pressure (CPR) preconditioning approach \cite{Wallis,Wallis85,wang}, the domain decomposition methods \cite{Liu2018,Sko13,haijian-cmame,haijian-jcp19}, and the algebraic or geometric multigrid algorithms \cite{bui17,bui18,qiao,wang}, are flexible and capable of addressing the inherent ill-conditioning of a complicated physics system, and hence have received increasingly more attentions in recent years. The focus of this paper is on the domain decomposition approach by which we propose a family of one-level or multilevel restricted additive Schwarz preconditioners. The original overlapping additive Schwarz (AS) method was introduced for the solution of  symmetric positive definite elliptic finite element problems, and was later extended to many other nonlinear or linear systems  \cite{cs99,Toselli}.  In the AS method, it follows the divide-and-conquer technique by recursively breaking down a problem into  more sub-problems of the same or related type, and its communication only occurs between neighboring subdomains during the restriction and extension processes.  Moreover, the additive Schwarz preconditioner does not require any splitting of the nonlinear system. Hence,  it can serve as a basis of efficient algorithms that  are naturally suitable for massively parallel computing.  In particular,  the approach can be combined with certain variants of restriction operators in a robust and efficient way, such as the so-called restricted additive Schwarz (RAS) method proposed by Xiao-Chuan Cai et al. in \cite{cs03,cs99}. Hence, it can significantly improve the physical feasibility of the solver and  substantial reduce the total computing time, and has proven to be very efficient in a variety of applications \cite{kong,Prudencio,Sko13,wang,haijian-cmame,haijian-jcp}. Recently, the overlapping restricted additive Schwarz method in conjunction with a variant of inexact Newton methods has been applied successfully to the two-phase flow  problem \cite{haijian-cmame} and the black oil model \cite{wang,haijian-jcp19}. It was demonstrated that the one-level Schwarz preconditioner is scalable in terms of the numbers of nonlinear and linear iterations, but not in terms of the total compute time {when the number of processors becomes large}, which means the requirement of the family of multilevel methods \cite{kong,Prudencio,sbg96,Toselli}.

In this work, to take advantage of modern supercomputers for large-scale reservoir simulations, we propose and develop the  overlapping restricted additive Schwarz preconditioner into a general multilevel framework for solving discrete systems coming from fully implicit discretizations. We show that,  with this new feature, our proposed approach is more robust and efficient for problems with highly heterogeneous media, and can scale optimally to a number of processors on the Tianhe-2 supercomputer. We would like to pointing out that designing a good strategy for the multilevel approach is both time-consuming and challenging, as it requires extensive knowledge of the general-purpose framework of interests, such as the selection of  coarse-to-fine mesh ratios, the restriction and interpolation operators, and the solvers for the smoothing and coarse-grid correction. To the best of our knowledge, very limited research has been conducted to apply the multilevel restricted additive Schwarz preconditioning technique for the fully implicit solution of petroleum reservoirs.

The outline of the paper is as follows.  In Section \ref{Model}, we introduce the governing equations of subsurface flows with Peng-Robinson equation of state, followed by the corresponding fully implicit discretization. In Section \ref{two-RAS}, a family of multilevel restricted additive Schwarz preconditioners, as the most important part of the fully implicit  solver, is presented in detail.  We show numerical results for  some 2D and 3D realistic benchmark problems in Section \ref{Numerical-scheme} to demonstrate the robustness and scalability of the proposed preconditioner. The paper is concluded in Section \ref{conclusion}.

\section{Mathematical model and discretizations}\label{Model}
The problem of interest in this study is the  compressible Darcy flow in a porous medium for the description of gas reservoir simulations \cite{Chen,Chen07,Firoozabadi}. Let $\Omega\in \mathbb{R}^d$ $(d =1,2,3)$ be the computational domain with the boundary $\Gamma=\partial \Omega$, and $t_f$ denote the final simulated time. The mass balance equation with the Darcy law for the real gas fluid  is defined by
\begin{equation}\label{mass conservation}
\left\{
\begin{array}{ll}
  \displaystyle \frac{\partial(\phi\rho)}{\partial{t}}+\nabla \cdot (\rho \mathbf{v})= q,  ~\mbox{in~} \Omega_{T}=\Omega\times (0, t_f], \\
 
   \displaystyle \mathbf{v}=-\frac{\mathbf{K}}{\mu}(\nabla p -\rho g\nabla h), \\
\end{array}
   \right.
  \end{equation}
where  $\phi$ is the porosity of the porous medium,   $\mathbf{K}$ is the permeability tensor,  $\mathbf{v}$ is the Darcy's velocity, $\rho$ denotes the density, $\mu$ is the viscosity of the flow, $p$ is the pressure, and $q$ is the external mass flow rate. The Peng-Robinson (PR) equation of state (EOS) \cite{Firoozabadi} is used to describe the density as a function of the composition, temperature and pressure:
\begin{equation}\label{equ:density}
\left\{
\begin{array}{ll}
\rho=cW,\\
c=\displaystyle\frac{p}{ZRT},\\
\end{array}
 \right.
\end{equation}
where $W$ is the molecular weight, $T$ is the temperature, $R$ is the gas constant, and $Z$  is the compressibility  factor of gas defined by
\begin{equation}\label{equ:PR}
Z^3-(1-B)Z^2+(A-3B^2-2B)Z-(AB-B^2-B^3)=0.
\end{equation}
Moreover,  $A$ and $B$ are the PR-EOS parameters defined as follows:
 \begin{equation*}
 \displaystyle A=\frac{ap}{R^2T^2}, ~ B=\frac{bp}{RT},
\end{equation*}
where the parameters $a$ and $b$ are modeled by imposing the criticality conditions
\begin{equation}\label{attrfactor}
\left\{
\begin{array}{ll}
\displaystyle a=a(T)=0.45724\cdot\frac{R^2T_c^2}{p_c}\left(1+m\left(1-\sqrt{T_r}\right)\right)^2,\\

\displaystyle b=0.0778\cdot\frac{RT_c}{p_c},
\end{array}
 \right.
\end{equation}
with  $T_c$ ($p_c$) being the specific pressure (temperature) in the critical state of the gas, and $T_r$ being the reduced temperature defined as ${T}/{T_c}$. In addition, the parameter $m$  in \eqref{attrfactor} is a fitting formula of the acentric factor $w$ of the substance:
\begin{equation}
 m= \left\{
   \begin{aligned}
 0.37464+1.54226w-0.26992w^2  \quad       0\leq w\leq 0.491,\\
 0.3796+1.485w-0.1644w^2+0.01667w^3 \quad  0.491\leq w\leq 2.0.\\
   \end{aligned}
   \right.
  \end{equation}
Here, the acentric factor $w$ is calculated by the following formula:
\begin{equation*}
 w= \frac{-ln(\frac{p_c}{1\;atm})-f^{(0)}(T_{br})}{f^{(1)}(T_{br})},
  \end{equation*}
where $f^0(T_{br})$ and $f^1(T_{br})$ are given by
\begin{equation*}
\left\{
\begin{aligned}
f^0(T_{br})&=\frac{-5.97616(1-T_{br})+1.29874(1-T_{br})^{\frac{3}{2}}-0.60394(1-T_{br})^{\frac{5}{2}}-1.06841(1-T_{br})^5}{T_{br}},\\
f^1(T_{br})&=\frac{-5.03365(1-T_{br})+1.11505(1-T_{br})^{\frac{ 3}{2}}-5.41217(1-T_{br})^{\frac{5}{2}}- 7.46628(1-T_{br})^5}{T_{br}},
\end{aligned}
\right.
\end{equation*}
with the normal boiling point
temperature $T_b$ and the reduced normal boiling point temperature defined as
$T_{br}={T_b}/{T_c}$.

In this study, we take the pressure $p$ as the primary variable.  Suppose that the porosity $\phi$ is not changed with the time, then the equations of the mathematical model can be rewritten  by substituting \eqref{equ:density} and \eqref{equ:PR} into \eqref{mass conservation} as follows:
\begin{equation}\label{gover_sys}
\frac{\phi p}{Z}\left(\frac{1}{p}-\frac{1}{Z}\frac{\partial Z}{\partial p}\right)\frac{\partial p}{\partial t}- \nabla \cdot \left(\frac{\mathbf{K}p}{\mu Z} \left(\nabla p -\frac{pW}{ZRT} g\nabla h\right)\right)- \frac{qRT}{W}=0,
  \end{equation}
where the relations are given by
\begin{equation*}
\left\{
   \begin{aligned}
&\frac{\partial Z}{\partial p}=-\frac{b}{RTC}Z^2-
    \left(\frac{ a}{R^2T^2C}-\frac{2b+6Bb}{RTC}\right)Z+\left(\frac{2 abp}{R^3T^3C}-\frac{2Bb+3B^2b}{RTC}\right),\\
    &C=3Z^2-2(1-B)Z +(A-2B-3B^2),\\
&Z^3-(1-B)Z^2+(A-2B-3B^2)Z-(AB-B^2-B^3)=0,
   \end{aligned}
   \right.
  \end{equation*}
with  the  initial condition  $p=p_0$.
Suppose the boundary of the computational domain $\Omega$ is composed of two parts $\partial\Omega=\Gamma^{D}+\Gamma^{N}$ with $\Gamma^{D}\cap\Gamma^{N}=\emptyset$.
The boundary conditions associated to the model problem \eqref{gover_sys} are
\begin{equation*}\label{boundary}
\left\{\begin{array}{ll}
p= p^{D} &\mbox{on}\, \Gamma^D,\\
\mathbf{v}\cdot\mathbf{n}=q^{N} &\mbox{on}\, \Gamma^N,\\
\end{array}\right.
\end{equation*}
where $\textbf{n}$ is the outward normal of the boundary $\partial\Omega$.  We remark that the  compressibility factor $Z$ as an intermediate variable is obtained by solving the algebraic cubic equation \eqref{equ:PR} with the primary variable $p$, see the references \cite{Chen,Chen07} for  the computation of $Z$ in details.

We employ a cell-centered finite difference (CCFD) method for the spatial discretization, for which the details  can be found in \cite{Chen,Chen07,Monteagudo07,haijian-siam}, and then a fully implicit scheme is applied for the time integration \cite{bui17,kong,wang,haijian-cmame}. For a given time-stepping sequence $0=t^{(1)}<t^{(2)}<...$, define the time step size $\Delta t^{(l)}=t^{(l+1)}-t^{(l)}$ and use superscript $(l)$ to denote the discretized evaluation at time level $t=t^{(l)}$. After spatially discretizing \eqref{gover_sys} by the CCFD scheme, we have a semi-discrete system
\begin{equation}\label{semi-discretized-system}
\frac{\phi p}{Z}\left(\frac{1}{p}-\frac{1}{Z}\frac{\partial Z}{\partial p}\right)\frac{\partial p}{\partial t}+\mathcal{F}\left(p^{(l+1)}\right)=0,
\end{equation}
where $\mathcal{F}\left(p^{(l+1)}\right)$ is the  operator of the spatial discretization at time level $t=t^{(l+1)}$ by using the CCFD method.  For the purpose of comparison,  we implement both the first-order backward Euler scheme (BDF-1) and  second-order backward differentiation formula (BDF-2) for the temporal integration of \eqref{semi-discretized-system}. For the BDF-1 scheme, the fully discretized system reads
\begin{equation}
\label{implicit-method-bdf1}
\left[\frac{\phi p}{Z}\left(\frac{1}{p}-\frac{1}{Z}\frac{\partial Z}{\partial p}\right)\right]^{(l+1)} \frac{p^{(l+1)}-p^{(l)}}{\triangle t^{(l)}}+\mathcal{F}\left(p^{(l+1)}\right)=0,
\end{equation}
where $p^{(l)}$ is the evaluation of $p$ at the $l^{th}$ time step with a  time step size $\triangle t^{(l)}$. And the fully discretized system for the BDF-2 method is
\begin{equation}
\label{implicit-method-bdf2}
\left[\frac{\phi p}{Z}\left(\frac{1}{p}-\frac{1}{Z}\frac{\partial Z}{\partial p}\right)\right]^{(l+1)} \frac{3p^{(l+1)}-4p^{(l)}+p^{(l-1)}}{2\triangle t^{(l)}}+\mathcal{F}\left(p^{(l+1)}\right)=0.
\end{equation}
Note that,  for the BDF-2 scheme, the BDF-1 method is used at the initial time step.

\begin{remark}
The employed CCFD discretization for \eqref{gover_sys} on rectangular meshes can be viewed as the mixed finite element method with Raviart-Thomas basis functions of lowest order equipped with the trapezoidal quadrature rule \cite{Chen,Chen07,Monteagudo07}. We also would like to pointing out that the proposed solver technologies introduced in the following section can be extended to solve the systems arising from the other spatial discretization schemes.
\end{remark}

\section{Multilevel restricted Schwarz preconditioners}\label{two-RAS}
When the fully implicit method is applied, a system of nonlinear algebraic equations after the temporal and spatial discretizations,
\begin{equation}
\label{equ:nonlinear-system}
  F(X)=0,
\end{equation}
is constructed and solved at each time step. In \eqref{equ:nonlinear-system},  the vector $F(X): \mathbb{R}^n\rightarrow \mathbb{R}^n$ is a given nonlinear vector-valued function arising from the residuals function in \eqref{implicit-method-bdf1} or  \eqref{implicit-method-bdf2} given by $F(X)=(F_1, F_2, \cdots, F_n)^{T}$ with $F_i=F_i(X_1, X_2,\cdots, X_n)^{T}$, and $X=(X_1, X_2,\cdots, X_n)^{T}$. In this study, the nonlinear system \eqref{equ:nonlinear-system}  is  solved by a parallel Newton-Krylov method with the family of domain decomposition type preconditioners \cite{Cai98,cs02a,fim5,Knoll04}.

Let the initial guess for Newton iterations $X^{0}\in\mathbb{R}^n$ at current time step be the solution of the previous time step, then the new solution $X^{k+1}$ is updated by the current solution $X^{k}$ and a Newton correction vector $s^{k}$ as follows:
\begin{equation}
  X^{k+1}=X^{k}+\lambda^{k}s^{k}, \quad k=0,1,2,\cdots,
  \label{newton_iterate}
\end{equation}
{
where $\lambda^{k}$ is the step length calculated by a cubic line search method \cite{ds} to satisfy
 \begin{equation*}\label{eq:line-search}
\begin{array}{ll}
\left\| F(X^k+\lambda^{k} s^{k})\right\|\leq (1-\alpha\lambda^{k})\left\|F(X^k)\right\|,
\end{array}
\end{equation*}
with  the parameter $\alpha$ being fixed to $10^{-4}$ in practice}. And  the Newton correction $s^{k}$  is obtained by solving the following  Jacobian linear system
\begin{equation}\label{jacobian1}
  J_k s^k =-F(X^k),
\end{equation}
with a Krylov subspace iteration method \cite{Sa03}.  Here,  $J_k=\nabla F(X^k)$ is the Jacobian matrix obtained from the current solution $X^{k}$. Since the corresponding linear system from the Newton iteration is large sparse and ill-conditioned,  the family of one-level or multilevel restricted Schwarz preconditioners is taken into account to accelerate the convergence of the linear system, i.e., we solve the following right-preconditioned Jacobian system
\begin{equation}\label{jacobian}
  J_k M_k^{-1} M_k s^k =-F(X^k),
\end{equation}
with an overlapping Schwarz preconditioned Krylov subspace method. In practice,  the Generalized Minimal RESidual (GMRES) algorithm as a Krylov subspace method is used for the linear solver.

 The accuracy of approximation \eqref{jacobian} is controlled by the relative and absolute tolerances $\eta^r_k$ and $\eta^a_k$ until the following convergence criterion is satisfied
\begin{equation}\label{equ:linear_stop}
  \| J_k s^k +F(X^k)\| \leq  \min\{\eta^r_k \| F(X^k)\|,\eta^a_k\}.
\end{equation}
We set the stopping conditions for the Newton iteration as
\begin{equation}\label{equ:nonlinear_stop}
||{F}(X^{k+1})|| \leq \max\{ \varepsilon_{r}
||{F}(X^{0})||, \varepsilon_{a} \},
\end{equation}
with the relative $\varepsilon_{r}$ and  absolute  $\varepsilon_{a}$ solver tolerances, respectively.

The most important component of a robust and scalable solver for solving the nonlinear or linear system is the choice of suitable  preconditioners.  In large-scale parallel computing, the additive Schwarz (AS) preconditioner, as a family of domain decomposition methods,  can help to improve the convergence and meanwhile is beneficial to the scalability of the linear solver,  which is the main focus of this paper.

\subsection{One-level Schwarz preconditioners}
Let $\Omega\subset \mathbb{R}^d$, $d=1,2,3$,  be the computational domain on which the PDE system \eqref{gover_sys} are defined. Then a spatial discretization is performed with a mesh $\Omega_h$ of the characteristic size $h>0$. To define the one-level Schwarz preconditioner,  we first divide $\Omega_h$ into nonoverlapping subdomains $\Omega_i$, $i=1,\cdots,N_p$, and then expand each $\Omega_i$ to $\Omega_i^{\delta}$, i.e., $\Omega_i\subset\Omega_i^{\delta}$, to obtain the overlapping partition. Here, the parameter  $\delta>0$  is the size of the overlap  defined as the minimum distance between $\partial \Omega_i^{\delta}$ and $\partial \Omega_i$, in the interior of $\Omega$. For boundary subdomains we simply cut off the part outside $\Omega$.  We also  denote  $H>0$ as the characteristic diameter of $\Omega_i$, as shown in the left panel of Figure \ref{fig:partition-domain}.

\definecolor{blackred}{RGB}{144,0,0}
\definecolor{mygreen}{RGB}{107,142,50}
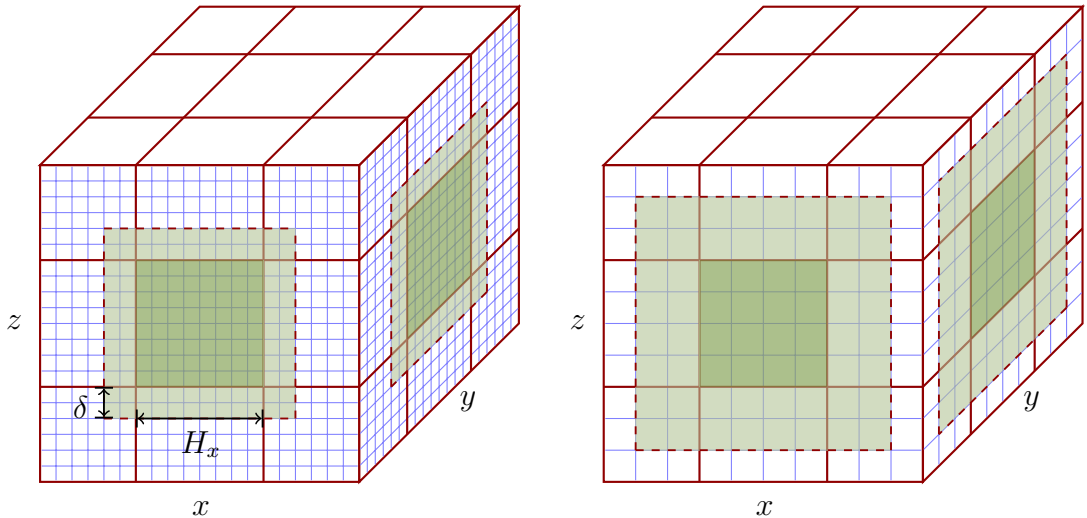
\begin{figure}[ht]
{\begin{center}
%%%%%%%%%%%%%%%%%%%%%%%%%%%%%%%%%%%%%%%%%%%%%%%%
\begin{tikzpicture}[x =1.0cm,y =1.0cm,z  = -0.5cm,
%x  = {(-0.75cm,-0.25cm)},
                    %y  = {(0.9659cm,-0.35882cm)},
                  %  y  = {(0.9659cm,-0.35882cm)},
                    %z  = {(0cm,1cm)},
                    scale = 1.05]
\tikzset{facestyle/.style={fill=lightgray,draw=black,very thin,line join=round}}

\begin{scope}[canvas is xy plane at z=4]
       %\path[facestyle] (0,0) rectangle (4,6);
       \fill[white,fill opacity=.9] (0,0) rectangle (4,4);
        \draw[blackred,thick] (0,0) rectangle (4,4);
        \draw[step=2.0mm, blue!50] (0.02,0.02) grid (3.98,3.98);
        %x-direction
         \draw[blackred,thick](1.2,0)node[right]{$\mathsf{~ }$} to (1.2,4);
         \draw[blackred,thick](2.8,0)node[right]{$\mathsf{~ }$} to (2.8,4);	
         %y-direction
         \draw[blackred,thick](0,1.2)node[right]{$\mathsf{~ }$} to (4,1.2);
         \draw[blackred,thick](0,2.8)node[right]{$\mathsf{~ }$} to (4,2.8);

        \fill[mygreen!80,opacity=.7] (1.2,1.2) rectangle (2.8,2.8);
         %\fill[gray!50,opacity=.7] (0.8,0.8) rectangle (3.2,1.2);
         \fill[mygreen!50,opacity=.6] (0.8,0.8) rectangle (3.2,1.2);
         \fill[mygreen!50,opacity=.6] (0.8,2.8) rectangle (3.2,3.2);
         \fill[mygreen!50,opacity=.6] (0.8,1.2) rectangle (1.2,2.8);
         \fill[mygreen!50,opacity=.6] (2.8,1.2) rectangle (3.2,2.8);

         \draw[blackred,line width =0.7pt,dashed](0.8,0.8)node[right]{$\mathsf{~ }$} to (0.8,3.2);
         \draw[blackred,line width =0.7pt,dashed](0.8,0.8)node[right]{$\mathsf{~ }$} to (3.2,0.8);
         \draw[blackred,line width =0.7pt,dashed](3.2,0.8)node[right]{$\mathsf{~ }$} to (3.2,3.2);
         \draw[blackred,line width =0.7pt,dashed](0.8,3.2)node[right]{$\mathsf{~ }$} to (3.2,3.2);

         \draw[|<->|,line width =0.7pt] (1.2,0.8) -- (2.8,0.8);
         \put(0,-46){\makebox(1.0,1.0){$H_x$}}
         \draw[|<->|,line width =0.7pt] (0.8,0.8) -- (0.8,1.2);
         \put(-45,-32){\makebox(1.0,1.0){$\delta$}}

\end{scope}
%% face  "right"
\begin{scope}[canvas is yz plane at x=4]
        \fill[white,fill opacity=.9] (0,0) rectangle (4,4);
        \draw[blackred,thick] (0,0) rectangle (4,4);
        \draw[step=2.0mm, blue!50] (0.02,0.02) grid (3.98,3.98);

         \draw[blackred,thick](1.2,0)node[right]{$\mathsf{~ }$} to (1.2,4);
         \draw[blackred,thick](2.8,0)node[right]{$\mathsf{~ }$} to (2.8,4);	
         %y-direction
         \draw[blackred,thick](0,1.2)node[right]{$\mathsf{~ }$} to (4,1.2);
         \draw[blackred,thick](0,2.8)node[right]{$\mathsf{~ }$} to (4,2.8);

        \fill[mygreen!80,opacity=.7] (1.2,1.2) rectangle (2.8,2.8);
         %\fill[gray!50,opacity=.7] (0.8,0.8) rectangle (3.2,1.2);
         \fill[mygreen!50,opacity=.6] (0.8,0.8) rectangle (3.2,1.2);
         \fill[mygreen!50,opacity=.6] (0.8,2.8) rectangle (3.2,3.2);
         \fill[mygreen!50,opacity=.6] (0.8,1.2) rectangle (1.2,2.8);
         \fill[mygreen!50,opacity=.6] (2.8,1.2) rectangle (3.2,2.8);

         \draw[blackred,line width =0.7pt,dashed](0.8,0.8)node[right]{$\mathsf{~ }$} to (0.8,3.2);
         \draw[blackred,line width =0.7pt,dashed](0.8,0.8)node[right]{$\mathsf{~ }$} to (3.2,0.8);
         \draw[blackred,line width =0.7pt,dashed](3.2,0.8)node[right]{$\mathsf{~ }$} to (3.2,3.2);
         \draw[blackred,line width =0.7pt,dashed](0.8,3.2)node[right]{$\mathsf{~ }$} to (3.2,3.2);

\end{scope}

           \put(0,-70){\makebox(1.0,1.0){$x$}}
           \put(-70,0){\makebox(1.0,1.0){$z$}}
           \put(100,-30){\makebox(1.0,1.0){$y$}}

% face "up"
\begin{scope}[canvas is xz plane at y=4]
       %\path[facestyle] (0,0) rectangle (4,6);
       \fill[white,fill opacity=.9] (0,0) rectangle (4,4);
        \draw[blackred,thick] (0,0) rectangle (4,4);
        %\draw[step=5mm, black] (0,0) grid (5.5,4);
                %x-direction
         \draw[blackred,thick](0,1.2)node[right]{$\mathsf{~ }$} to (4,1.2);	
         \draw[blackred,thick](0,2.8)node[right]{$\mathsf{~ }$} to (4,2.8);	
         %y-direction
         \draw[blackred,thick](1.2,0)node[right]{$\mathsf{~ }$} to (1.2,4);	
         \draw[blackred,thick](2.8,0)node[right]{$\mathsf{~ }$} to (2.8,4);
%         \draw[blue,thick,dashed](4.25,0)node[right]{$\mathsf{~ }$} to (4.25,4);
\end{scope}
\end{tikzpicture}
\hspace{0.4cm}
%%%%%%%%%%%%%%%%%%%%%%%%%%%%%%%%%%%%%%%%%%%
\begin{tikzpicture}[x =1.0cm,y =1.0cm,z  = -0.5cm,
%x  = {(-0.75cm,-0.25cm)},
                    %y  = {(0.9659cm,-0.35882cm)},
                  %  y  = {(0.9659cm,-0.35882cm)},
                    %z  = {(0cm,1cm)},
                    scale = 1.05]
% style of faces
\tikzset{facestyle/.style={fill=lightgray,draw=black,very thin,line join=round}}

\begin{scope}[canvas is xy plane at z=4]

      \fill[white,fill opacity=.9] (0,0) rectangle (4,4);
        \draw[blackred,thick] (0,0) rectangle (4,4);
        \draw[step=4.0mm, blue!50] (0.02,0.02) grid (3.98,3.98);
        %x-direction
         \draw[blackred,thick](1.2,0)node[right]{$\mathsf{~ }$} to (1.2,4); 	
         \draw[blackred,thick](2.8,0)node[right]{$\mathsf{~ }$} to (2.8,4);	
         %y-direction
         \draw[blackred,thick](0,1.2)node[right]{$\mathsf{~ }$} to (4,1.2);
         \draw[blackred,thick](0,2.8)node[right]{$\mathsf{~ }$} to (4,2.8);

        \fill[mygreen!80,opacity=.7] (1.2,1.2) rectangle (2.8,2.8);
         \fill[mygreen!50,opacity=.6] (0.4,0.4) rectangle (3.6,1.2);
         \fill[mygreen!50,opacity=.6] (0.4,2.8) rectangle (3.6,3.6);
         \fill[mygreen!50,opacity=.6] (0.4,1.2) rectangle (1.2,2.8);
         \fill[mygreen!50,opacity=.6] (2.8,1.2) rectangle (3.6,2.8);

         \draw[blackred,line width =0.7pt,dashed](0.4,0.4)node[right]{$\mathsf{~ }$} to (0.4,3.6);
         \draw[blackred,line width =0.7pt,dashed](0.4,0.4)node[right]{$\mathsf{~ }$} to (3.6,0.4);
         \draw[blackred,line width =0.7pt,dashed](3.6,0.4)node[right]{$\mathsf{~ }$} to (3.6,3.6);
         \draw[blackred,line width =0.7pt,dashed](0.4,3.6)node[right]{$\mathsf{~ }$} to (3.6,3.6);

\end{scope}
%% face  "right"
\begin{scope}[canvas is yz plane at x=4]
         \fill[white,fill opacity=.9] (0,0) rectangle (4,4);
        \draw[blackred,thick] (0,0) rectangle (4,4);
        \draw[step=4.0mm, blue!50] (0.02,0.02) grid (3.98,3.98);

         \draw[blackred,thick](1.2,0)node[right]{$\mathsf{~ }$} to (1.2,4);
         \draw[blackred,thick](2.8,0)node[right]{$\mathsf{~ }$} to (2.8,4);	
         %y-direction
         \draw[blackred,thick](0,1.2)node[right]{$\mathsf{~ }$} to (4,1.2);
         \draw[blackred,thick](0,2.8)node[right]{$\mathsf{~ }$} to (4,2.8);

        \fill[mygreen!80,opacity=.7] (1.2,1.2) rectangle (2.8,2.8);
         \fill[mygreen!50,opacity=.6] (0.4,0.4) rectangle (3.6,1.2);
         \fill[mygreen!50,opacity=.6] (0.4,2.8) rectangle (3.6,3.6);
         \fill[mygreen!50,opacity=.6] (0.4,1.2) rectangle (1.2,2.8);
         \fill[mygreen!50,opacity=.6] (2.8,1.2) rectangle (3.6,2.8);

         \draw[blackred,line width =0.7pt,dashed](0.4,0.4)node[right]{$\mathsf{~ }$} to (0.4,3.6);
         \draw[blackred,line width =0.7pt,dashed](0.4,0.4)node[right]{$\mathsf{~ }$} to (3.6,0.4);
         \draw[blackred,line width =0.7pt,dashed](3.6,0.4)node[right]{$\mathsf{~ }$} to (3.6,3.6);
         \draw[blackred,line width =0.7pt,dashed](0.4,3.6)node[right]{$\mathsf{~ }$} to (3.6,3.6);

\end{scope}

% face "up"
\begin{scope}[canvas is xz plane at y=4]
       %\path[facestyle] (0,0) rectangle (4,6);
      \fill[white,fill opacity=.9] (0,0) rectangle (4,4);
        \draw[blackred,thick] (0,0) rectangle (4,4);
         \draw[blackred,thick](0,1.2)node[right]{$\mathsf{~ }$} to (4,1.2);	
         \draw[blackred,thick](0,2.8)node[right]{$\mathsf{~ }$} to (4,2.8);	
         \draw[blackred,thick](1.2,0)node[right]{$\mathsf{~ }$} to (1.2,4);	
         \draw[blackred,thick](2.8,0)node[right]{$\mathsf{~ }$} to (2.8,4);
\end{scope}

           \put(0,-70){\makebox(1.0,1.0){$x$}}
           \put(-70,0){\makebox(1.0,1.0){$z$}}
           \put(100,-30){\makebox(1.0,1.0){$y$}}

\end{tikzpicture}
\end{center}}
\caption{A demonstrated partition of  domain decompositions with geometry preserving coarse meshes. In the figure,  the red solid lines indicate the partition of the 3D domain into $3\times3\times3$ non-overlapping subdomains of size $H_x\times H_y\times H_z$, the dotted lines show the extended boundary of an overlapping subdomain with $\delta=2$. The left (right) panel of the figure denotes the partition with  a fine (coarse) mesh  size $h$ ($2h$), respectively.}
 \label{fig:partition-domain}
\end{figure}

Let $N$ and $N_i$ denote the number of degrees of freedom associated to $\Omega$ and $\Omega_j^{\delta}$, respectively. Let $A\in \mathbb{R}^{N \times N}$ be the  Jacobian matrix of the linear system defined on a mesh $\Omega_h$
 \begin{equation}\label{lin_system1}
AX=b.
\end{equation}
Then we can define the matrices ${R}^{\delta}_i$ and ${R}^0_i$ as the restriction operator from $\Omega_h$ to its overlapping and non-overlapping subdomains as follows: Its  element $({R}^{\delta}_i)_{l_1,l_2}$ is either (a) an identity block, if the integer indices
$1\leq l_1\leq N_i$ and $1\leq l_2\leq N$ are related to the same mesh point and this mesh point belongs to $\Omega_j^{\delta}$,
or (b) a zero block, otherwise. The multiplication of ${R}_i^{\delta}$ with a $N \times 1$ vector generates
a smaller $N_i\times1$ vector by discarding all elements corresponding to mesh points outside $\Omega_i^{\delta}$.
 The matrix ${R}_i^0\in \mathbb{R}^{N_i\times N}$ is defined in a similar way,  the only difference to the operator ${R}^{\delta}_i$ is that  its application to a $N\times 1$ vector also zeros out all of those elements corresponding to
mesh points on $\Omega_i^{\delta}\backslash \Omega_i$.
 Then the classical one-level additive Schwarz (AS, \cite{Toselli}) preconditioner is defined as
\begin{equation}
\begin{array}{ll}
{M}^{-1}_{\delta,\delta} =\sum\limits_{i=1}^{N_p}({R}^{\delta}_i)^T
{A}_i^{-1}{R}^\delta_i.
\end{array}
\label{equ:as-class}
\end{equation}
with ${A}_i= {R}^\delta_i {A}({R}^\delta_i)^T$ and $N_p$ is the number of subdomains, which is the same as the number of processors. In addition to that,  there are two modified approaches of the one-level additive Schwarz preconditioner
that may have some potential advantages for parallel computing. The first version is the left restricted additive Schwarz (left-RAS, \cite{cs99}) method defined by
\begin{equation}
\begin{array}{ll}
{M}^{-1}_{0,\delta} =\sum\limits_{i=1}^{N_p}({R}^0_i)^T
{A}_i^{-1}{R}^\delta_i.
\end{array}
\label{equ:as-left}
\end{equation}
and the other modification to the original method is the right restricted additive Schwarz (right-RAS, \cite{cs03}) preconditioner as follow:
\begin{equation}
\begin{array}{ll}
{M}^{-1}_{\delta,0} =\sum\limits_{i=1}^{N_p}({R}^\delta_i)^T
{A}_i^{-1}{R}^0_i.
\end{array}
\label{equ:as-right}
\end{equation}
In the above preconditioners, we use a sparse LU factorization or incomplete LU (ILU) factorization  method to solve the subdomain linear system corresponding to the matrix ${A}_i^{-1}$.   In the following, we will denote a Schwarz preconditioner simply by $M^{-1}_{h}$ defined on the mesh $\Omega_h$ with the characteristic size $h$, when the distinction is not important.

\subsection{Multilevel additive Schwarz preconditioners}
To improve the scalability  and robustness of the one-level additive Schwarz preconditioners, especially when a large number of processors is used,  we employ the family of multilevel Schwarz preconditioners.  Multilevel Schwarz preconditioners are obtained  by combining the single level preconditioner $M^{-1}_{h}$ assigned to each level, as shown in Figure \ref{fig:partition-domain}. For the description of multilevel Schwarz preconditioners \cite{kong,Prudencio,sbg96}, we use the index $j=0,1,...,L-1$ to designate any of the $L\geq 2$ levels.  The meshes from coarse to fine are denoted by $\Omega_{h_j}$, and the corresponding matrices and vectors are denoted $A_{h_j}$ and $X_{h_j}$.  Let us denote $I_{h_j}$ as the identity operator defined on the level $j$, and  the restriction operator from the level $j$ to the level $j-1$ be defined by
\begin{equation}
\mathcal{I}_j^{j-1}: \mathbb{R}^{N_j} \rightarrow \mathbb{R}^{N_{j-1}},
\label{equ:restriction}
\end{equation}
where $N_j$ and $N_{j-1}$ denote the number of degrees of freedom associated to $\Omega_{h_j}$ and $\Omega_{h_{j-1}}$.  Moreover,
\begin{equation}
\label{equ:interpolation}
\mathcal{I}_{j-1}^{j}: \mathbb{R}^{N_{j-1}} \rightarrow \mathbb{R}^{N_{j}}
\end{equation}
is the interpolation operator from the level $j-1$ to the level $j$.

For the convenience of introduction, we present the construction of  the proposed multilevel  Schwarz preconditioners from the view of the multigrid (MG) V-cycle algorithm \cite{kong,Prudencio,sbg96,Toselli}. More precisely speaking, in this sense, at each level $j>0$, the Schwarz preconditioned Richardson method works as the pre-smoother and  post-smoother,  i.e.,  $M^{-1}_{h_j}$  preconditioning the $\mu_i\geq0$ presmoother iterations and $M^{-1}_{h_j}$ preconditioning the $\nu_i\geq0$ postsmoother iterations. In the general multigrid V-cycle framework,   let $X^k$ be the current solution for the  linear system \eqref{lin_system1},  the new solution is computed by the iteration $X^{k+1}=MG(b,L,X^{k})$, as described in Algorithm \ref{vcyclealgorithm}.

\begin{algorithm}
\caption{Multigrid (MG) V-cycle algorithm.}\label{vcyclealgorithm}
\begin{algorithmic}[1]
\Procedure{$X_{h_j}=MG$}{$b_{h_j},j,X_{h_j}$}
\If{$j=0$}
\State Solve $A_{h_0}X_{h_0}=b_{h_0}$\Comment{Coarsest correction}
\Else
\State Smooth $\mu_i$ times $A_{h_j}X_{h_j}=b_{h_j}$:\Comment{Presmoothing}
\State $(b_{h_j}-A_{h_j}X_{h_j})=(I_{h_j}-A_{h_j}M^{-1}_{h_j})^{\mu_i}(b_{h_j}-A_{h_j}X_{h_j})$;
\State $b_{h_{j-1}}=\mathcal{I}_j^{j-1}(b_{h_j}-A_{h_j}X_{h_j})$; \Comment{Residual restriction}
\State  $X_{h_{j-1}}=MG(b_{h_{j-1}},j-1,0)$;\Comment{Recursivity}
\State $X_{h_j}=X_{h_j}+\mathcal{I}_{j-1}^{j} X_{h_{j-1}}$;\Comment{Correction interpolation}
\State Smooth $\nu_i$ times $A_{h_j}X_{h_j}=b_{h_j}$:\Comment{Postsmoothing}
\State $(b_{h_j}-A_{h_j}X_{h_j})=(I_{h_j}-A_{h_j}M^{-1}_{h_j})^{\nu_i}(b_{h_j}-A_{h_j}X_{h_j})$; \EndIf
\State \textbf{return} $X_{h_j}$
\EndProcedure
\end{algorithmic}
\end{algorithm}

\subsubsection{With application to the two-level case}
In the following, we restrict ourselves to the two-level case, and use the geometry preserving coarse mesh that shares the same boundary geometry with the fine mesh, as shown in  Figure \ref{fig:partition-domain}. The V-cycle two-level Schwarz preconditioner $M_{V}^{-1}$, i.e, $L=2$ in Algorithm \ref{vcyclealgorithm}, is constructed by combining the fine  level $M^{-1}_{h_1}$ and the coarse level $M^{-1}_{h_0}$ preconditioners  as
follows:
\begin{equation}\label{two_level_preconditioner}
M_{V}^{-1}=A_{h_1}^{-1}\left[I_{h_1}-(I_{h_1}-A_{h_1}M^{-1}_{h_1})^{\nu_1}(I_{h_1}-A_{h_1}\mathcal{I}_0^{1}M^{-1}_{h_0}\mathcal{I}_1^0)(I_{h_1}-A_{h_1}M^{-1}_{h_1})^{\mu_1}\right].
\end{equation}
More precisely, if the pre- and post-smoothing parameters are fixed to  $\mu_i=1$ and  $\nu_i=1$, then the matrix-vector product of the two-level Schwarz preconditioner  in \eqref{two_level_preconditioner} with any given vector $e_{h_1}=M_{V}^{-1}r_{h_1}$ is obtained by the following three steps:
\begin{equation}\label{equ:v-cycle}
\left\{
\begin{array}{ll}
\displaystyle{ e_{h_1}^{1/3}=M^{-1}_{h_1}r_{h_1} },\\
\\
\displaystyle{ e_{h_1}^{2/3}=e_{h_1}^{1/3}+\mathcal{I}^{1}_{0}M^{-1}_{h_0}\mathcal{I}^{0}_{1} (r_{h_1}-A_{h_1}e_{h_1}^{1/3}),}\\
\\
\displaystyle{ e_{h_1}=e_{h_1}^{2/3}+M^{-1}_{h_1} (r_{h_1}-A_{h_1}e_{h_1}^{2/3}).}
\end{array}
 \right.
\end{equation}

There are some modified versions of the V-cycle two-level Schwarz preconditioner that can be used for parallel computing.  If we set $\mu_1=0$ and $\nu_1=1$ in \eqref{two_level_preconditioner}, then we obtain the left-Kaskade Schwarz method defined by
\begin{equation}\label{equ:left-Kas}
M_{left-Kas}^{-1}=\mathcal{I}_0^{1}M^{-1}_{h_0}\mathcal{I}_1^0+M^{-1}_{h_1}-M^{-1}_{h_1}A_{h_1}\mathcal{I}_0^{1}M^{-1}_{h_0}\mathcal{I}_1^0.
\end{equation}
And the other modification to the original method is the right-Kaskade Schwarz preconditioner with $\mu_1=1$ and $\nu_1=0$ in \eqref{two_level_preconditioner} as follow:
\begin{equation}\label{equ:right-Kas}
M_{right-Kas}^{-1}=M^{-1}_{h_1}+\mathcal{I}_0^{1}M^{-1}_{h_0}\mathcal{I}_1^0-\mathcal{I}_0^{1}M^{-1}_{h_0}\mathcal{I}_1^0A_{h_1}M^{-1}_{h_1}.
\end{equation}
Moreover, we can define other two hybrid versions of two-level Schwarz algorithms. The first two-level method  is the pure additive two-level Schwarz preconditioner as
\begin{equation}\label{equ:two-additive}
M_{additve}^{-1}=\mathcal{I}_0^{1}M^{-1}_{h_0}\mathcal{I}_1^0+M^{-1}_{h_1},
\end{equation}
and the coarse-level  only  type two-level Schwarz preconditioner
\begin{equation}\label{equ:coarse-level}
M_{coarse}^{-1}=\mathcal{I}_0^{1}M^{-1}_{h_0}\mathcal{I}_1^0.
\end{equation}
The motivation of the above two-level preconditioners, including, additively or multiplicatively, coarse preconditioners to an existing fine mesh preconditioner, is to make the overall solver scalable with respect to the number of processors or the number of subdomains. Hence, some numerical results will be shown later to compare the performance of these two-level preconditioners.

%\begin{remark}\label{remark-1}
%In the traditional multigrid approach,  one can certainly apply some strategies to choose the pre- and post-smoothing times $\mu_i, \nu_i$. However,  in practice, the potential improvement in performance is usually offset by the cost to construct the approximation,  our experiments suggest that there is no benefit to include the second pre- and post-smoothing times in the fine level of the one-level preconditioning for the implicit solution of the model problem, and a similar performance can be observed in the references \cite{kong,Prudencio}.
%\end{remark}

\begin{remark}\label{remark-2}
In the multilevel Schwarz preconditioners,  the use of a coarse level helps the exchange of information, and the following linear problem defined on the coarsest mesh needs to be solved
\begin{equation}
A_{h_0}X_{h_0}=b_{h_0},
\label{equ:coarse-system}
\end{equation}
with the use of the one-level  additive Schwarz preconditioned GMRES, owing to the face that it is too large and expensive for some direct methods in the applications of large-scale simulations.  When  an iterative method is used for solving the  linear system on the coarsest level \eqref{equ:coarse-system}, the overall preconditioner is an iterative procedure, and  the preconditioner changes from iteration to iteration. Hence, when the multilevel Schwarz method is applied,  a flexible version of GMRES (fGMRES, \cite{Sa03}) is  used for the solution of the outer linear system \eqref{lin_system1}.
\end{remark}

\begin{remark}\label{remark-3}
When the classical  additive Schwarz preconditioner \eqref{equ:as-class} is applied to symmetric positive definite
systems arising from the discretization of elliptical problems, then the condition number $\kappa$ of the preconditioned linear system satisfies $\kappa\leq C(1+H/\delta)/{H^2}$ for the one-level method and $\kappa\leq C(1+H/\delta)$ for the two-level method, where the parameter  $C$ is independent of $\delta$, $H$, and $h$, see the references \cite{sbg96,Toselli}. However, these condition number estimates can neither be applied to the restricted additive Schwarz preconditioners, nor adapted to the family of  non-elliptic systems like our model problem, and thereby there are very little theoretical literatures on the convergence of multilevel restricted Schwarz preconditioners for this system. The scalability tests in this paper will provide more understanding of restricted type domain decomposition methods for the hyperbolic problem.
\end{remark}

\subsubsection{Selection of interpolation operators}
 Algorithm \ref{vcyclealgorithm} provides a general framework for choosing the coarse-to-fine mesh ratio, the restriction and interpolation operators, and the solvers for the smoothing and coarse-grid correction.
And in this among them, choosing the right restriction and interpolation operators at each level is very important for the overall performance of the preconditioner in terms of the trade-off between rate of convergence and the cost of each  iteration. Generally speaking, in most cases the overall solver makes a profit of the addition of coarse preconditioners with the decrease of the total number of linear iterations.  However, for some classes of important problems such as reservoir  simulation,  we found that the number of iterations does not decrease as expected,
if the operators are not chosen properly. After some attempts, we figure out the source of this phenomenon from the cell-centered spatial discretization. When the cell-centered scheme is involved, it  happens that a coarse mesh point
does not coincide with any fine mesh point, as shown in Figure \ref{two-level-domain}. Hence, in the rest of this paper, we focus on the some strategies for the selection of  interpolation operators, i.e, the first-order,  second-order and third-order  interpolation schemes.

\definecolor{pinkeyellow}{RGB}{255,227,132}
\definecolor{mycolor1}{RGB}{128,118,105}
\begin{figure*}[ht]
{\begin{center}
\begin{tikzpicture}[scale=.8,every node/.style={minimum size=1cm},on grid]		

\begin{scope}[ yshift=-155,every node/.append style={
            yslant=0.5,xslant=-1},yslant=0.5,xslant=-1
            ]
        \fill[white,fill opacity=.9] (0,0) rectangle (5,5);
        \draw[black,very thick] (0,0) rectangle (5,5);
        \draw[step=5mm, black] (0,0) grid (5,5);
        %\draw[step=1mm, red!50,thin] (3,1) grid (4,2);  %Nested Grid
        \fill[red] (1.00,1.00) rectangle (1.5,1.5);
        \fill[green] (1.5,1.5) rectangle (2.0,2.0);
        \fill[mycolor1] (1.5,1.0) rectangle (2.0,1.5);
        \fill[yellow] (1.0,1.5) rectangle (1.5,2.0);
        \fill[mygreen!70,opacity=.6] (0.5,0.5) rectangle (1.0,2.5);
        \fill[mygreen!70,opacity=.6] (1.0,0.5) rectangle (2.5,1.0);
        \fill[mygreen!70,opacity=.6] (2.0,1.0) rectangle (2.5,2.5);
        \fill[mygreen!70,opacity=.6] (1.0,2.0) rectangle (2.0,2.5);
\end{scope}

   \begin{scope}[scale=0.7, yshift=-240,xshift=273,every node/.append style={
    	    yslant=0,xslant=0},yslant=0,xslant=0
    	             ]
    \draw[step=10mm, black] (0,0) grid (4,4);
    \foreach \c in {(0,0), (1,0), (2,0),(3,0),(0,1),(1,1),(2,1),(3,1),(0,1),(1,1),(2,1),(3,1),(0,2),(1,2),(2,2),(3,2),(0,3),(1,3),(2,3),(3,3)}
    \fill \c + (0.5,0.5) circle (0.1);
    \end{scope}

\begin{scope}[            yshift=-83,every node/.append style={
            yslant=0.5,xslant=-1},yslant=0.5,xslant=-1
            ]
        % opacity to prevent graphical interference
        \fill[white,fill opacity=.9] (0,0) rectangle (5,5);
        \draw[black,very thick] (0,0) rectangle (5,5);
        \draw[step=5mm, black] (0,0) grid (5,5);
        %\draw[step=1mm, red!50,thin] (3,1) grid (4,2);  %Nested Grid
        \fill[red] (1.00,1.00) rectangle (1.5,1.5);
        \fill[green] (1.5,1.5) rectangle (2.0,2.0);
        \fill[mycolor1] (1.5,1.0) rectangle (2.0,1.5);
        \fill[yellow] (1.0,1.5) rectangle (1.5,2.0);
\end{scope}

   \begin{scope}[scale=0.7, yshift=-83,xshift=300,every node/.append style={
    	    yslant=0,xslant=0},yslant=0,xslant=0
    	             ]
    \draw[step=10mm, black] (0,0) grid (2,2);
    \foreach \c in {(0,0), (1,0),  (0,1),(1,1)}
    \fill \c + (0.5,0.5) circle (0.1);
    \end{scope}

        \begin{scope}[yshift=2,every node/.append style={
    	    yslant=0.5,xslant=-1},yslant=0.5,xslant=-1
    	             ]
        \fill[white,fill opacity=1.0] (0,0) rectangle (5,5);
        \draw[black,very thick] (0,0) rectangle (5,5);
        \draw[step=5mm, dashed,red!50,thin] (0,0) grid (1,1);
        \draw[step=10mm, black] (0,0) grid (5,5);
        \fill[blue] (1.00,1.00) rectangle (2,2);
    \end{scope}

\begin{scope}[scale=0.7, yshift=25,xshift=300,every node/.append style={
    	    yslant=0,xslant=0},yslant=0,xslant=0
    	             ]
     \draw [scale=2.0, black] (0,0) grid (1,1);
     \foreach \c in {(0,0)}
     \fill \c + (1.0,1.0) circle (0.1);
\end{scope}

     \draw[-latex,thick](7.0,1.3)node[right]{$\mathsf{~ }$}
        to[out=180,in=90] (0,1.5);
    \draw[-latex,thick](7.0,-1.4)node[right]{$\mathsf{~ }$}
        to[out=180,in=90] (0,-1.4);	
    \draw[-latex,thick](6.3,-4.5)node[right]{$\mathsf{~ }$}
        to[out=180,in=90] (0,-4);	
\end{tikzpicture}
\end{center}}
\caption{Partitions for a 2D spatial domain  with the coarse-to-fine mesh ratio 1:2. } \label{two-level-domain}
\end{figure*}
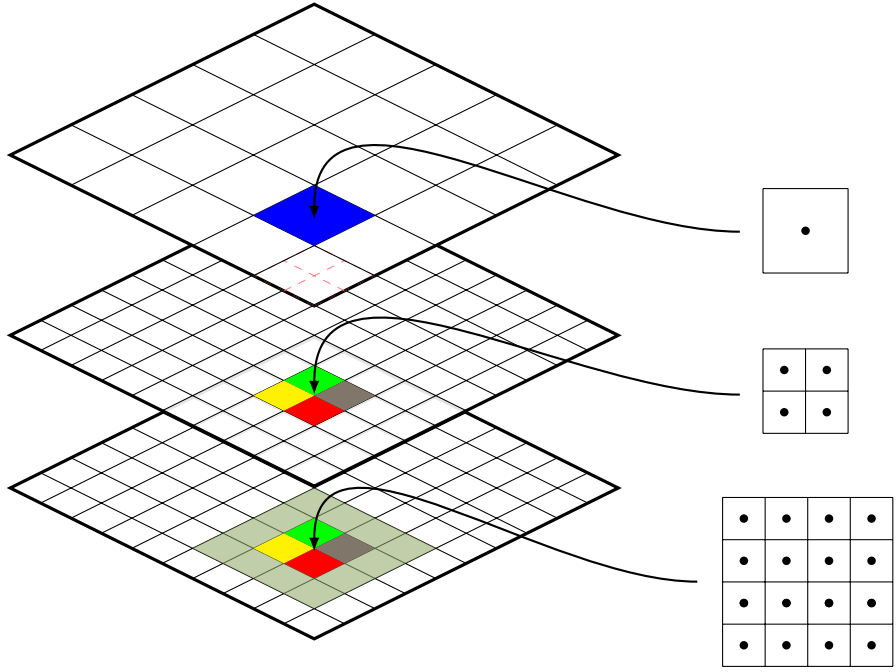

 Let  $\Omega$ be the computational domain covered with $N_x\times N_y$ mesh cells. Then we consider $(x, y)$ is a point at the position of the rectangular subdivision $[x_i,x_{i+1}]\times[y_j,y_{j+1}]$ with $i=1, \dots, N_x$,  $j=1, \dots, N_y$, and $f(x,y)$ is the interpolation point value on this point.  The nearest neighbor interpolation operator, which  approximate the interpolation point data of the nearest node according to the shortest distance between the interpolation point and the sample point in space, is a first-order method defined by
 \begin{equation}\label{equ:inter-frist}
f(x,y)= \sum_{l=0}^1\sum_{m=0}^1f\left(x_{i+l},y_{j+m}\right)W_{i+l,j+m},
\end{equation}
where $f(x_{j+l},y_{k+m})$ is the interpolation point value located at the  point $(x_{j+l},y_{k+m})$, and $W_{i+l,j+m}$ is a weight defined as
\begin{equation*}
W_{i+l,j+m}=\left\{\begin{aligned}
1, \quad &\mbox{if} ~\|(x,y)-(x_{i+l},y_{j+m})\|~\mbox{is  minimum}, \\
0, \quad &\mbox{others}.
\end{aligned}
\right.
\end{equation*}
The bilinear interpolation operator, which  utilizes the weighted average of the nearest neighboring values $f(x_{i+l},y_{j+m})$ to approximately generate a interpolation point value. Due to the point $(x, y)$ being in the subdivision $[x_i,x_{i+1}]\times[y_j,y_{j+1}]$,  is a second-order interpolation technique defined as follows:
\begin{equation}
\label{equ:inter-second}
f(x,y)=\sum_{l=0}^1\sum_{m=0}^1f\left(x_{i+l},y_{j+m}\right)W_{i+l}W_{j+m},
\end{equation}
where the weight for the bilinear interpolation scheme is defined as:
\begin{equation*}
\left\{
\begin{aligned}
W_{i+l}&=\left(\frac{2x-x_i-x_{i+1}}{x_{i+1}-x_{i}}\right)l+\frac{x_{i+1}-x}{x_{i+1}-x_{i}},\\
W_{j+m}&=\left(\frac{2y-y_j-y_{j+1}}{y_{j+1}-y_{j}}\right)m+\frac{y_{j+1}-y}{y_{j+1}-y_{j}}.
\end{aligned}
\right.
\end{equation*}
 In contrast to the bilinear interpolation, which only takes 4 ponits $(2\times2)$ into account, the bicubic interpolation use 16 points $(4\times4)$ and is a third-order interpolation scheme defined as follows:
 \begin{equation}
 \label{equ:inter-third}
   f(x,y)=\sum_{l=-1}^2\sum_{m=-1}^2f\left(x_{i+l},y_{j+m}\right)W\left(\frac{x-x_{i+l}}{x_{i+1}-x_{i}}\right)W\left(\frac{y-y_{j+m}}{y_{j+1}-y_{j}}\right),
 \end{equation}
where the weight for the cubic convolution interpolation is defined as
\begin{equation*}
W(s)=\left\{\begin{aligned}
&\frac{3}{2}|s|^3-\frac{5}{2}|s|^2+1,\quad &0< |s|< 1, \\
&-\frac{1}{2}|s|^3+\frac{5}{2}|s|^2-4|s|+2, \quad&1<|s|< 2,  \\
&\;0,  \quad &2<|s|.
\end{aligned}
\right.
\end{equation*}
We remark that the two dimensional computational  domain is used in the above description only for the ease of demonstration.

%Then, we focus on constructing the restriction operator with second and third-order accuracies. Let $(x, y)$ be a point in the rectangular subdivision $[x_i,x_{i+1}]\times[y_j,y_{j+1}]$ divided from the computational domain $\Omega$, and $f(x,y)$ be the interpolation point value on this point.  Then the two-dimensional restriction function with the  second-order accuracy is
%\begin{equation*}
%f(x,y)=\sum_{l=0}^1\sum_{m=0}^1f(x_{i+l},y_{j+m})W_{i+l}W_{j+m},
%\end{equation*}
%where the weight for the this restriction function is defined as
%\begin{equation*}
%\left\{
%\begin{aligned}
%&W_{i}=\frac{1}{2},W_{i+1}=\frac{1}{2},\\
%&W_{j}=\frac{1}{2},W_{j+1}=\frac{1}{2}.
%\end{aligned}
%\right.
%\end{equation*}
%Moreover, the restriction function with the third-order accuracy takes the following form
% \begin{equation*}
%   f(x,y)=\sum_{l=-1}^2\sum_{m=-1}^2f(x_{i+l},y_{j+m})W_{i+l}W_{j+m},
% \end{equation*}
%where the weight for the restriction function is defined as
%\begin{equation*}
%\left\{
%\begin{aligned}
%&W_{i-1}=-\frac{1}{16},W_{i}=\frac{9}{16},W_{i+1}=\frac{9}{16},W_{i+2}=-\frac{1}{16}\\
%&W_{j-1}=-\frac{1}{16},W_{j}=\frac{9}{16},W_{j+1}=\frac{9}{16},W_{j+2}=-\frac{1}{16}.
%\end{aligned}
%\right.
%\end{equation*}

%Here $s=\frac{x-x_{i+l}}{x_{i+1}-x_{i}}$ or $\frac{y-y_{j+m}}{y_{j+1}-y_{j}}$ and the boundary condition are given third-order approximation.

\section{Numerical experiments}\label{Numerical-scheme}
In this section, we implement the proposed algorithm described in the previous sections  using the open-source Portable, Extensible Toolkit for Scientific computation (PETSc) \cite{petsc}, which is built on the top of Message Passing Interface (MPI), and investigate the numerical behavior and parallel performance of the newly proposed fully implicit solver with a variety of test cases. 

%(c) the comparison of the performance of the multilevel restricted Schwarz preconditioners with respect to different parameters; and  (d) the parallel performance of the fully implicit solver.

%In the fully implicit solver, we set the tolerances as follows. For the Newton iteration, an absolute (relative) tolerance of $10^{-10}$ ($10^{-6}$) is utilized. For the linear iteration, the linear systems are solved by the Schwarz preconditioned GMRES method with absolute and relative tolerances of $10^{-8}$ and $10^{-5}$, respectively.

\subsection{\red{Robustness and efficiency of the solver}}
\red{The focus of the subsection is on the robustness and efficiency of the proposed algorithm for both standard benchmarks and realistic problems in highly heterogeneous media.}  Unless otherwise specified, the values of physical parameters used in the test cases are set as follows: $\phi=0.2$, $\mu=11.067\times10^{-3}  \; \mbox{cp}$, $R= 8.3147295\; \mbox{J}/(\mbox{mol}\cdot \mbox{K})$, $W=16\times10^{-3} \; \mbox{kg/mol}$, $T=298 \; \mbox{K}$, $p_c= 4.604\times10^6 \; \mbox{pa}$, $T_c=190.58 \; \mbox{K}$, and $T_b=111.67 \; \mbox{K}$.

\begin{figure}[H]
\centering
    \subfloat[Isotropic case]{\includegraphics[width=0.41\linewidth]{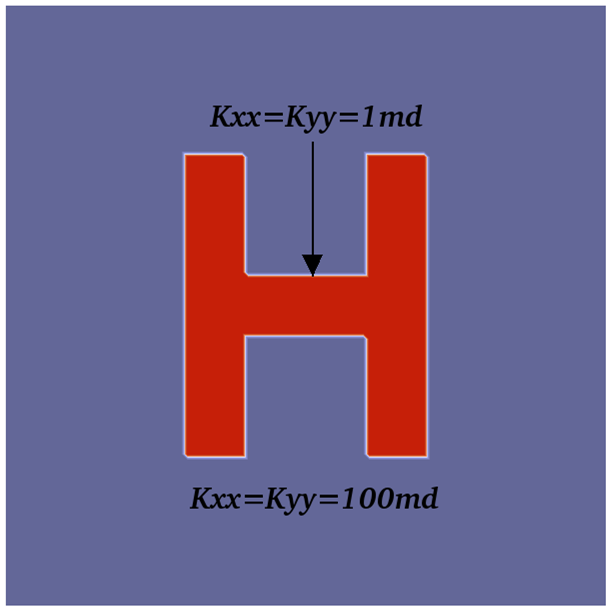}}
    \subfloat[Anisotropic case]{\includegraphics[width=0.41\linewidth]{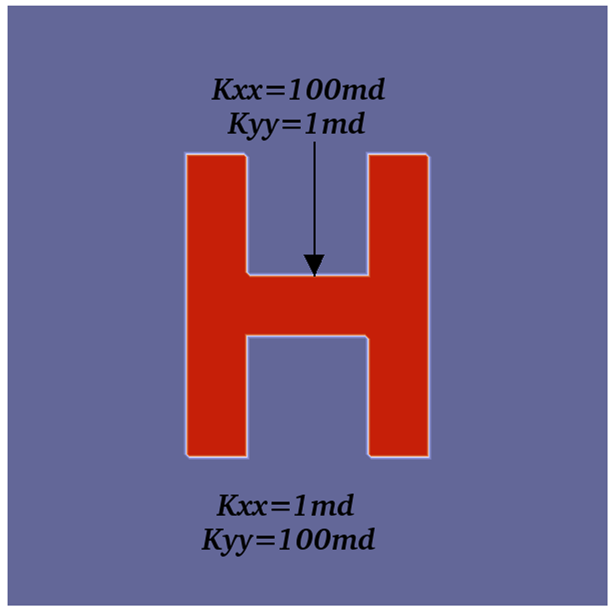}}    
      
    \subfloat[Isotropic case]{\includegraphics[width=0.41\linewidth]{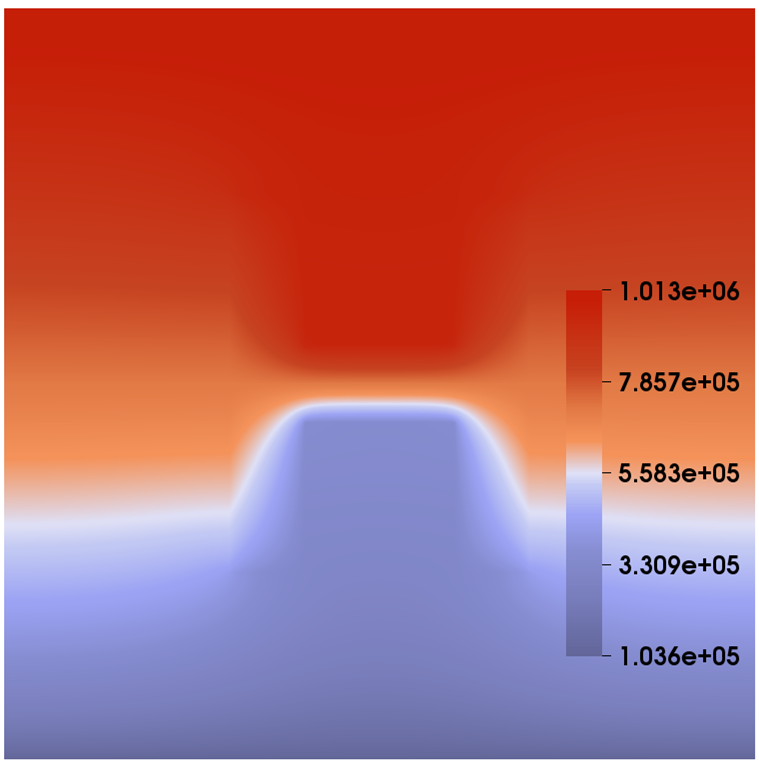}}
    \subfloat[Anisotropic case]{\includegraphics[width=0.411\linewidth]{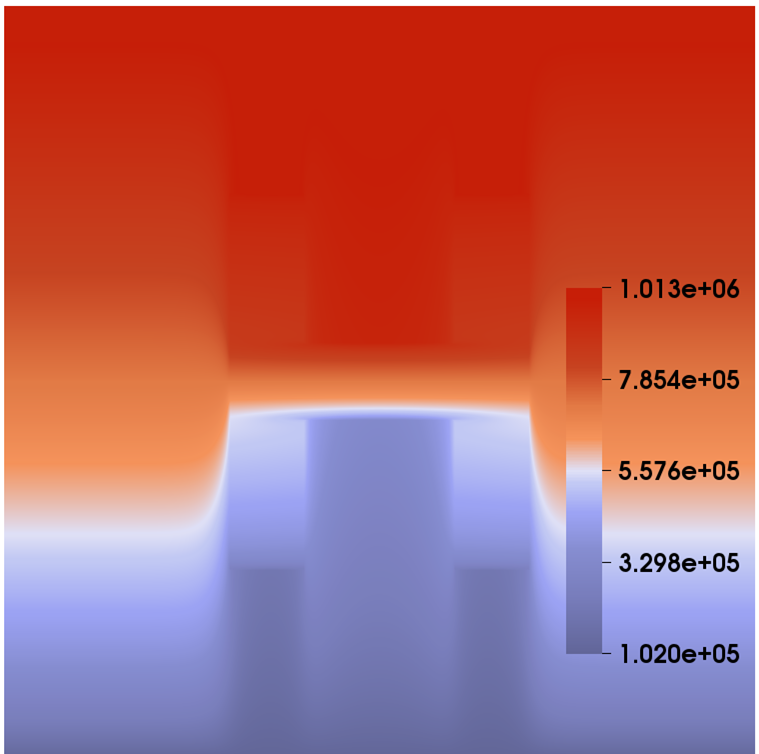}}
\caption{\red{Contour plots of the solution under the  heterogeneous isotropic and anisotropic mediums  for Case-1.}}
\label{case2_perm1}
\end{figure}

\red{We first present results from a 2D test case (denoted as Case-1) in a horizontal layer with the permeability tensor $\mathbf{K}$ in \eqref{mass conservation} as follows:
\begin{equation*}
\mathbf{K}=\left[\begin{array}{cc}
 k_{xx} & 0 \\
 0 & k_{yy}\\
\end{array}
 \right].
\label{2.4}
\end{equation*}
 With the media being horizontal, we neglect the effect of the gravity.  In the  configuration, the distribution of permeabilities includes two domains with different isotropy,  as shown in Figure \ref{case2_perm1}.  The computational  domain is  100 meters long and 100 meters wide.} In the simulation, we assume that the left and right boundaries of the domain are impermeable. Then we flood the system by gas from the top to the bottom, i.e., we set $p_{in}=10$ atm at the top boundary and set $p_{out}=1$ atm at the right boundary. And there is no injection/extraction inside the domain. Compared with the previous example,  in this test case  there is an $H$-shape zone and the value of permeabilities has a huge jump inside and outside the zone, which brings about even greater challenges to the fully implicit solver. In the test, the simulation is performed on a $512\times512$ mesh,  the  time step size is fixed to  $\Delta t=0.1$ day,  and the simulation is stopped at $10$ day. {Figure \ref{case2_perm1} also illustrates the contour plots of the pressure.} Table \ref{Heterogeneous_isotropic} shows  the performance of the proposed fully  implicit solver with respect to different permeability configurations. It is clearly seen that the simulation spends more computing time when the isotropic medium is used, which attributes to the increase in nonlinearity of the problem that affects the number of linear iterations.

%\begin{figure}[H]
%\centering
%    \subfloat[Isotropic case]{\includegraphics[width=0.45\linewidth]{./fig/exam3p.png}}
%    \hspace{0.02cm}
%    \subfloat[Anisotropic case]{\includegraphics[width=0.45\linewidth]{./fig/exam4p.png}}\\
%     \subfloat[$x=35.2$ m]{\includegraphics[width=0.45\linewidth]{./fig/case2x=352m}}
%    \subfloat[$x=50$ m]{\includegraphics[width=0.45\linewidth]{./fig/case2x=50m}}
%\caption{Contour plots of the pressure  and its profile curves along $x=35.2\; \mbox{m}$ and $x=50\; \mbox{m}$ for Case-1.  }
%\label{case3_pressure}
%\end{figure}

\begin{table}[hbp]
\begin{center}
 \caption{{\red{A comparison of the fully implicit solver with different permeability configurations for Case-1.}}  }
\begin{tabular}{|l|c|c|}
\hline
Medium type &{Isotropic case}&{Anisotropic case
}\\
\hline
%\multirow{3}{*}{Case-2}&Average  nonlinear iterations &2.2  &2.2\\
%&Average  linear iterations & 58.1  &45.7\\
%&Execution time (second)  &180.1  &160.0\\
%\hline
Average   nonlinear iterations &2.2  &2.2\\
Average  linear iterations& 54.7  &41.5\\
Execution time (second)  & 172.0  &151.1\\
\hline
\end{tabular}
\label{Heterogeneous_isotropic}
\end{center}
\end{table}

\begin{figure}[H]
\centering
    \includegraphics[width=0.4\linewidth]{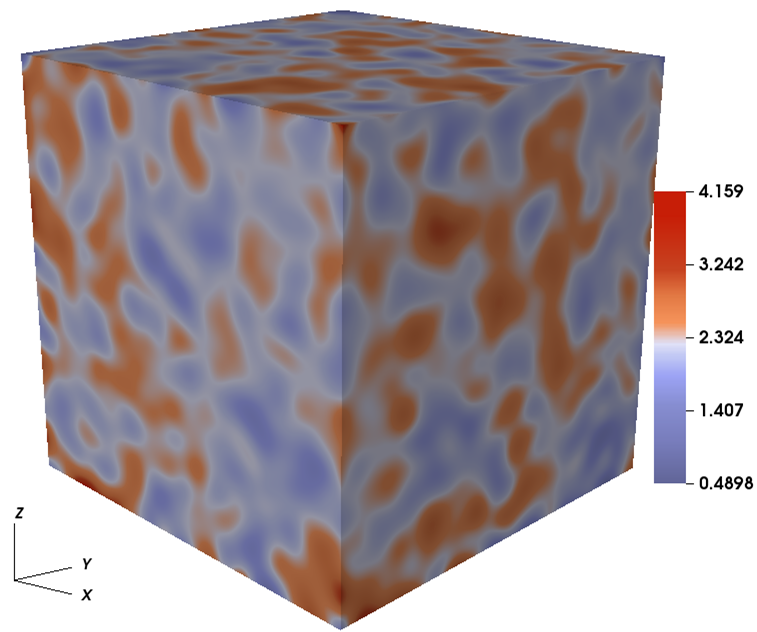}
\caption{Permeability field for Case-2. In the figure, we use the logarithmic scale for the contour plot of the permeability. }
\label{case3_perm}
\end{figure}

%\subsubsection{3D test cases with PR-EOS}
\red{In the following, the experiment  is conducted to simulate some 3D test cases. The focus is  on the flow model with the medium being highly heterogeneous and isotropic, which significantly increases the nonlinearity of the system and imposes a severe challenge on the fully implicit  solver. We first consider a 3D domain with dimension $\Omega=(0,100\; \mbox{m})^3$, in which the permeabilities of the porous medium are random, denoted as Case-2.} The random distribution of permeability with the range  [3.1,14426.2] is generated by a geostatistical model using the open source toolbox MRST \cite{mrst}, as shown in Figure \ref{case3_perm}. {In the test, two different flow configurations are used for the injection boundary, i.e., the fixed injection pressure and the fixed injection rate. In the fixed pressure configuration, the flow with the pressure of 10 atm is injected at the part of the left-hand side $x=0\; \mbox{m},y\in[0,10\; \mbox{m}],z\in[0,10\; \mbox{m}]$;  while,  for the case of the fixed injection rate, the flow is injected from the same sub-boundary with the Darcy's velocity $\mathbf{v}=5\times10^{-4}~\mbox{m/s}$.  The fluid flows from the right-hand side  of the domain with a fixed pressure $p_{out}=1$ atm, and no-flow boundary conditions are imposed on the other boundaries of the domain.}  The initial pressure of $p_0=1$ atm is specified in the whole domain and the parameters required for this example are consistent with the previous examples.

\begin{figure}[htp]
\centering
    \subfloat[t=1 day]{\includegraphics[width=0.33\linewidth]{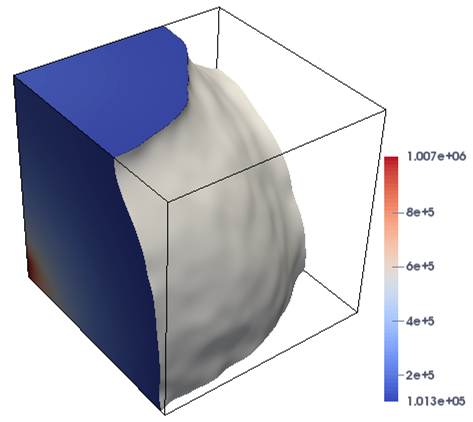}}
    \subfloat[t=2 day]{\includegraphics[width=0.33\linewidth]{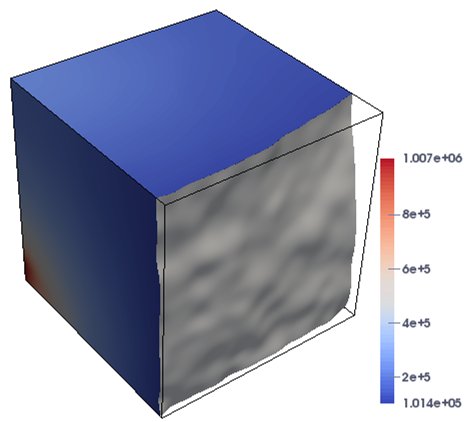}}
    \subfloat[t=5 day]{\includegraphics[width=0.33\linewidth]{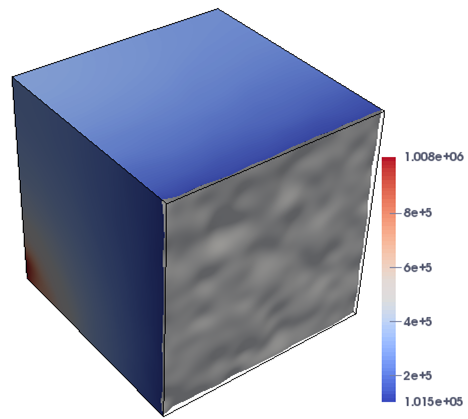}}\\
    \subfloat[t=1 day]{\includegraphics[width=0.33\linewidth]{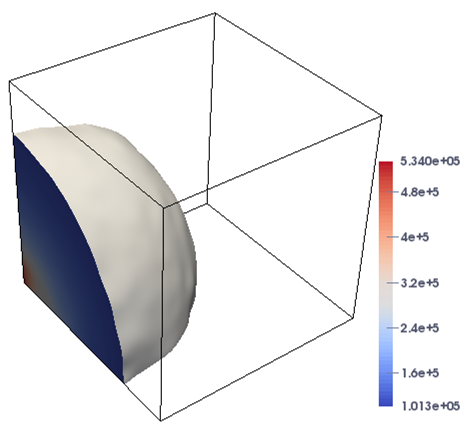}}
    \subfloat[t=2 day]{\includegraphics[width=0.33\linewidth]{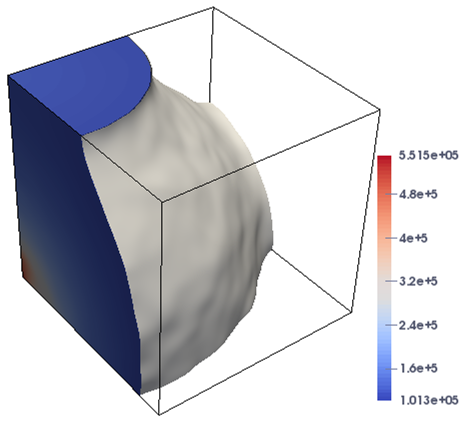}}
    \subfloat[t=5 day]{\includegraphics[width=0.33\linewidth]{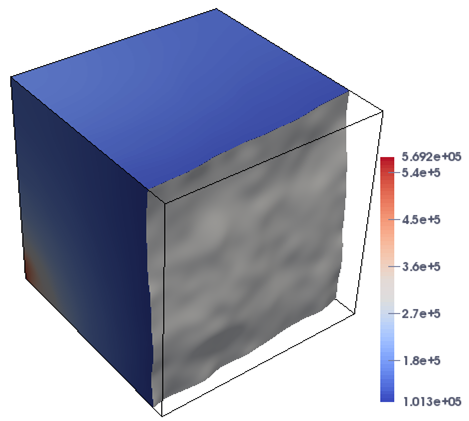}}
%    \subfloat[t=1 day]{\includegraphics[width=0.35\linewidth]{./fig/vel0_0008exam3_t=11.png}}
%    \subfloat[t=2 day]{\includegraphics[width=0.35\linewidth]{./fig/vel0_0008exam3_t=21.png}}
%    \subfloat[t=5 day]{\includegraphics[width=0.35\linewidth]{./fig/vel0_0008exam3_t=51.png}}
\caption{{A comparison of the pressure  profiles  at different times  for Case-2 with the fixed pressure (first row)  or velocity (second row) boundary.} }
\label{3D_pressure}
\end{figure}

In this case, the domain is highly heterogeneous,  leading to the increase in nonlinearity of the problem that is challenging for the numerical techniques. {Figure \ref{3D_pressure} shows the plots of the pressure profiles  for Case-2 at different times under the two flow configurations. It is demonstrated from the figures that the proposed approach successfully resolves the different stages of the simulation, and the flow is close to reach the breakthrough at the time $t=5$ days.} In Table  \ref{random_permeaiblity}, we present the history of the values of nonlinear and linear iterations and the total execution time for the proposed solver at different times. The simulation is carried out on a $128\times128\times128$ mesh and the time step size is $\Delta t=0.1$ day. It can be observed from Table  \ref{random_permeaiblity} that the number of nonlinear and linear iterations has barely changed, while the  total execution time increases constantly with the advance of time as expected, which displays the robustness and effectiveness of the implicit solver to handle the highly variety of physical  parameters in the model problem.

\begin{table}[hbp]
\begin{center}
{
 \caption{{A comparison of the fully implicit solver at different times for Case-2 with the fixed pressure or velocity boundary.} }
\begin{tabular}{ |c |l| c| c| c| c| }
\hline
Boundary Type&End of the simulation (day) &{$t=1$ }&{$t=2$}&{$t=3$ }&{ t=5 }\\
\hline
\multirow{3}{*}{Fixed pressure}&Average nonlinear iterations &3.0  & 2.9 &2.9 &2.9\\
&Average linear iterations & 22.6  &23.8 &24.3&24.5 \\
&Execution time (second)  &645.3  &1095.9&1545.2&2455.6\\
\cline{1-6}
\multirow{3}{*}{Fixed velocity}&Average nonlinear iterations &2.8&2.7&2.7&2.6\\
&Average linear iterations &24.8&26.3&27.1&27.9\\
&Execution time (second)  &519.5&965.4&1325.1&2142.9\\
\hline
%Fixed velocity&Average nonlinear iterations &2.9&2.8&2.8&2.6\\
%$u=0.0008m/s$&Average linear iterations &25.9&27.2&28.0&28.9\\
%&Execution time (second)  &563.6&1033.8&1459.6&2290.4\\
%\hline
\end{tabular}
\label{random_permeaiblity}
}
\end{center}
\end{table}

\begin{figure}[hbp]
\centering
    \subfloat[$log_{10}$-permeability field]{\includegraphics[width=0.45\linewidth]{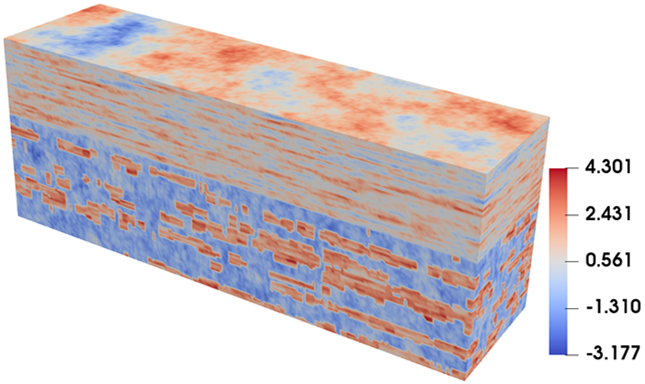}}
    \hspace{0.02cm}
    \subfloat[Porosity field]{\includegraphics[width=0.45\linewidth]{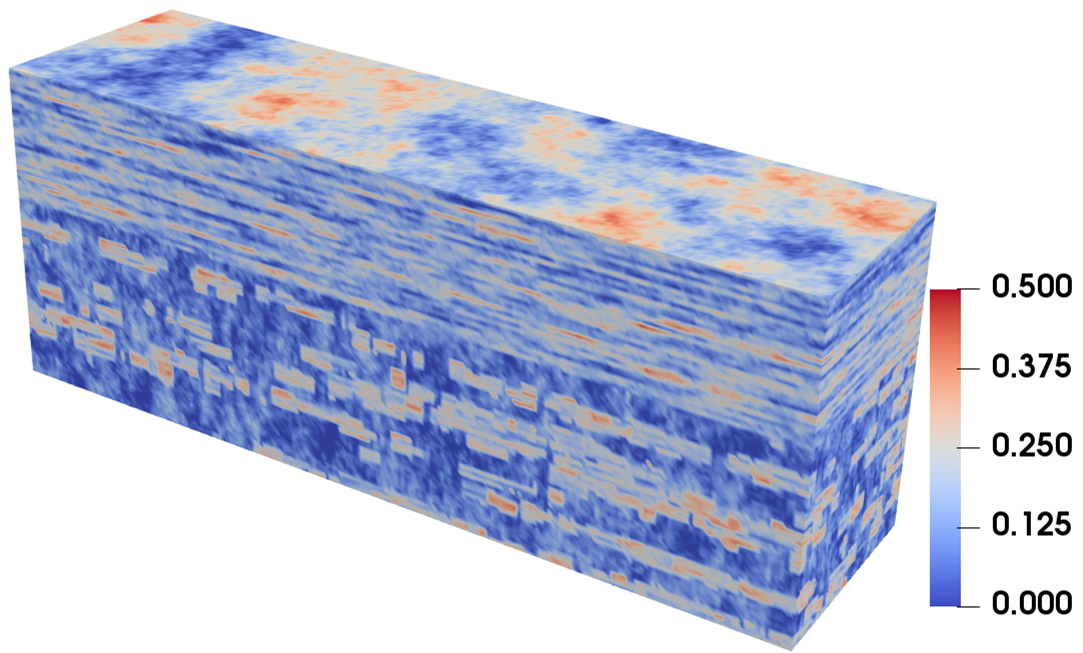}}\\
\caption{ {Permeability and porosity fields  for the SPE10 benchmark.}  }
\label{spe10permporo}
\end{figure}

In the all of the above test cases, the porosity $\phi$ in the model problem \eqref{mass conservation} is assumed to a constant.   In the next 3D test case (Case-3) of this  subsection,  we import the porosity and the permeability  from the Tenth SPE Comparative Project (SPE10) \cite{spe10} as an example of a realistic realization with geological and petrophysical properties, in which the porosity and the permeability are capable of variation with the change of the position. It is a classic and challenging benchmark problem for reservoir due to  highly heterogeneous permeabilities and porosities.  As shown in Figure \ref{spe10permporo}, the permeability is characterized by variations of more than six orders of magnitude and is ranged from $6.65\times10^{-4}$ to $2\times10^4$, and the porosity scale ranges from 0 to 0.5. The 3D domain dimensions are 1200 ft long $\times$ 2200 ft wide $\times$ 170 ft thick. In the test, the boundary condition for pressure on the left-hand side of the domain is uniformly imposed to $p_{in}=10$ atm and on the opposite side is imposed to $p_{out}=1$ atm. No-flow boundary condition is set to other boundary of the domain. Other used parameters are the same as the previous case.

\begin{figure}[hbp]
\centering
    \subfloat[t=1 day]{\includegraphics[width=0.45\linewidth]{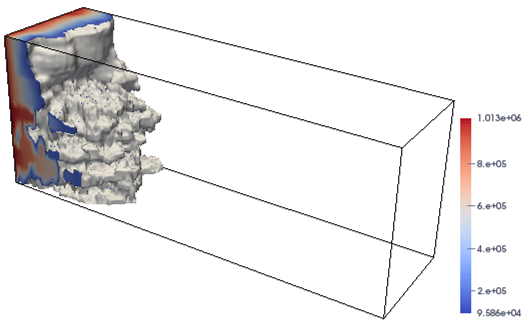}}
    \subfloat[t=3 day]{\includegraphics[width=0.45\linewidth]{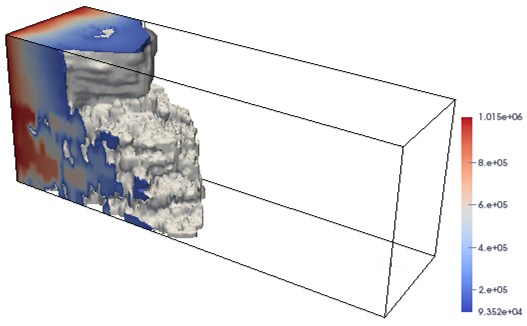}}\\
    \subfloat[t=6 day]{\includegraphics[width=0.45\linewidth]{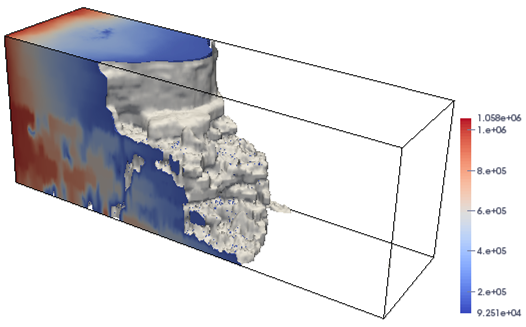}}
    \subfloat[t=10 day]{\includegraphics[width=0.45\linewidth]{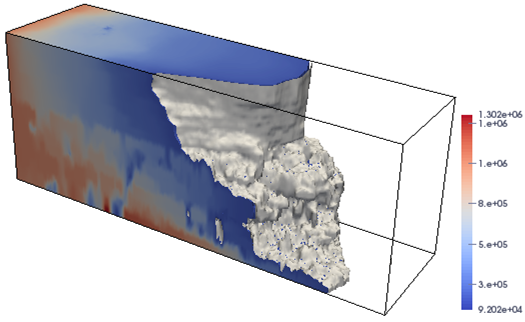}}
\caption{{A comparison of the pressure  profiles  at different times for the SPE10 benchmark.}}
\label{spe10_pressure}
\end{figure}

{In Figure \ref{spe10_pressure}, we display the distributions for the pressure  when the simulation is finished at different time points $t=1,3,6,10$ days.} In the simulation, the mesh is $60\times220\times85$ and the time step size is fixed to  $0.1$ day.  The results shown in Figure \ref{spe10_pressure} demonstrate that  the proposed approach successfully resolves the rapid and abrupt evolution of the simulation at different stages, while keeping the solver in a robust  and efficient way. Finally, we analyze the behavior of the proposed fully-implicit  method when the time step size $\Delta t$ is changed. In the test,  we again run the SPE-10 model on a fixed $60\times220\times85$ mesh. The simulation is stopped at $t=3$ year. The results on the average numbers of nonlinear  and linear iterations as well as the total compute time  are summarized in Table \ref{spe10_timesteps}.  The results in the table clearly indicate that the combination of nonlinear and linear iterations works well for even very large value of time steps. The implicit approach converges for all time steps and is unconditionally stable. In addition to that, we also notice that, as the time step size $\Delta t$ decreases, the average number of nonlinear and linear iterations become smaller, whereas the total computing time increases. This behavior is somehow expected for the fully implicit approach in a variety of applications \cite{bui17,kong,wang,haijian-cmame20,haijian-cmame}.

\begin{table}[H]
\begin{center}
 \caption{The effect of different time step sizes in the fully implicit solver for SPE10. }
\begin{tabular}{|l|c|c|c|c|c|}
\hline
Time step size $\Delta t$ &{0.06}&{0.075}&{0.1}&{0.2} &{0.3}\\
\hline
Number of time steps &50 &40 &30 &15 &10   \\
Average nonlinear iterations &3.5 &3.6&3.8 &4.4& 4.5\\
Average linear iterations & 45.4 & 49.3& 52.1& 67.5 &80.0\\
Execution time (second)  &1424.2  &1191.7&951.4 &579.8 &408.8\\
\hline
\end{tabular}
\label{spe10_timesteps}
\end{center}
\end{table}

\subsection{Performance of Schwarz preconditioners}  \label{sec:Schwraz}
In this subsection, we focus on the parallel performance of the proposed fully implicit solver with respect to the one-level or multilevel Schwarz preconditioner by Case-1, {Case-2}, and a new 3D test case. In the new 3D case, {denoted as Case-4},  the computational domain is $\Omega=(0,100 \; \mbox{m})^3$ with  the medium being heterogeneous and isotropic, and the distribution of permeability in each layer $z\in[25\; \mbox{m},75\; \mbox{m}]$ is the same as the distribution of Case-1. We flood the system by gas from the behind face ($x\in[0,100\; \mbox{m}]$, $y=100\; \mbox{m}$, $z\in[0,100\; \mbox{m}]$) to the front face, i.e., we set $p_{in}=10$ atm at the behind face and set $p_{out}=1$ atm at the front face. And there is no injection/extraction inside the domain.

In the test, we use the following stopping parameters and notations, unless specified otherwise. The relative  and absolute tolerances in \eqref{equ:nonlinear_stop} for Newton iterations are set to  $10^{-6}$ and $10^{-10}$, respectively. The linear systems are solved by the one-level or two-level Schwarz preconditioned GMRES method with absolute and relative tolerances of $10^{-8}$ and $10^{-5}$ in  \eqref{equ:linear_stop}, except for the coarse solve of the two-level preconditioner, for which we use $2\times 10^{-1}$ as the relative stopping tolerance. The restart value of the GMRES method  is fixed as 30.  In the tables of the following tests, the symbol ``$N_p$" denotes the number of processor cores that is the same as the number of subdomains, ``N. It" stands for the average number of Newton iterations per time step, ``L. It"  the average number of the one-level or two-level Schwarz preconditioned GMRES iterations per Newton iteration, and ``Time"  the total computing time in seconds.

\subsubsection{One-level Schwarz preconditioners}
Under the framework of one-level Schwarz preconditioners,  we study the performance of the fully implicit solver,
when different types of the additive Schwarz preconditioners are employed and  several important parameters such as the subdomain solver and the overlapping size  are taken into consideration.

We first look at  the influence of  subdomain solvers. In this test, we focus on the classical AS preconditioner  \eqref{equ:as-class} and fix the overlapping factor to $\delta=1$. The ILU factorization with different levels of fill-in and the full LU factorization are considered for the subdomain solvers. The simulation of Case-1 is applied to a fixed $512\times512$ mesh using $N_p=16$ processors. The time step size is 0.1 day, and we run the 2D case for the first 20 time steps. The simulation of Case-4 is applied to a fixed $128\times128\times128$ mesh by using $N_p=64$ processors.  The time step size is 0.1 day, and we run this case for the first 10 time steps. The numerical results are summarized in Table \ref{different-subpreconditioner}. We observe from the table that the number of nonlinear iterations is not sensitive to the choice of subdomain solvers.  The linear system converges with less iterations when the LU factorization method is used as the subdomain solver; while the ILU approach is beneficial to the AS preconditioner in terms of the computing time, especially for the 3D test case.

\begin{table}[H]
\begin{center}
 \caption{The impact of different subdomain solvers in the one-level Schwarz preconditioner.}
\begin{tabular}{|c|ccc|ccc|}
\hline
Case&\multicolumn{3}{c|}{Case-1} &\multicolumn{3}{c|}{Case-4}\\
\cline{1-4}\cline{5-7}
Solver & N. It &L. It &Time & N. It &L. It &Time\\
\hline
  \multirow{1}{*}{LU}
        &3.0   & 52.9  	&37.1     &4.2  &32.9    &954.8   \\
  \multirow{1}{*}{ILU(0)}
        &3.0   &99.1 	&34.3     &4.3  & 161.7   &282.4  \\
      \multirow{1}{*}{ILU(1)}
        &3.0  &75.7  	&29.9    &4.3 &  98.0  & 253.0 \\
      \multirow{1}{*}{ILU(2)}
        &3.0   &73.9 	&32.1    &4.3  &76.1  &239.6  \\
       \multirow{1}{*}{ILU(4)}
   	&3.0   & 66.2  &33.8    &4.3  &51.0  & 245.8 \\
\hline
\end{tabular}
\label{different-subpreconditioner}
\end{center}
\end{table}

{Moreover, in Table \ref{tab:ilu-case4}, we also compare the impact of different subdomain solvers for Case-2, when the fixed pressure and velocity  boundaries are imposed for the flow configurations. In the test, the simulation is carried out on a $128\times128\times128$ mesh and the time step size is fixed to $\Delta t=0.1$ day by using $N_p=64$ processors.  It is clearly illustrated from the results that the choice of ILU(2) is still optimal for the 3D flow problem under the random permeability case with different boundary conditions, when the flow undergoes different stages in the complex simulation as shown in Figure \ref{3D_pressure}.}

\begin{table}[htp]
\begin{center}
{
{\footnotesize
 \caption{{The impact of different subdomain solvers in the one-level Schwarz preconditioner for Case-2 with  the fixed pressure boundary (denoted as ``Model-A") and the fixed velocity  boundary (denoted as ``Model-B").}}
\begin{tabular}{|c|ccc|ccc|ccc|ccc|}
\hline
Solver&\multicolumn{3}{c|}{t=1} &\multicolumn{3}{c|}{t=2} &\multicolumn{3}{c|}{t=3} &\multicolumn{3}{c|}{t=5}\\
\cline{1-4}\cline{5-7}\cline{8-10}\cline{11-13}
Model-A & N. It &L. It &Time & N. It &L. It &Time  & N. It &L. It &Time  & N. It &L. It &Time \\
\hline
  \multirow{1}{*}{LU}
        &3.0&33.0&685.2&2.9&36.6& 1167.0  &2.9&38.7&1652.7&2.9&41.2&2625.9  \\
  \multirow{1}{*}{ILU(0)}
        &3.0&161.5&435.3&2.9&183.0& 760.1  &2.9&195.3&1086.8&2.9&208.3&1742.7   \\
      \multirow{1}{*}{ILU(1)}
        &3.0&100.6&415.5&2.9&110.3& 715.9   &2.9&116.4& 1015.4&2.9&123.6& 1616.7 \\
      \multirow{1}{*}{ILU(2)}
       &3.0&80.0&403.4&2.9& 88.9& 695.0 &2.9&93.3&991.1&2.9&98.4&1572.5   \\
       \multirow{1}{*}{ILU(4)}
   	&3.0&52.4&417.8&2.9&58.8&725.1  &2.9&62.1&1030.5&2.9&65.7&1638.9 \\
\hline
Model-B & N. It &L. It &Time & N. It &L. It &Time  & N. It &L. It &Time  & N. It &L. It &Time \\
\hline
  \multirow{1}{*}{LU}
        &2.8&39.3&563.7&2.7&42.9&1046.9 &2.7&45.2&1438.8&2.6& 47.7&2335.0\\
  \multirow{1}{*}{ILU(0)}
        &2.8&178.5&366.0&2.7&190.4&679.2 &2.7&201.8&934.4&2.6&214.9&1518.6\\
      \multirow{1}{*}{ILU(1)}
       &2.8&113.9&349.6&2.7&124.0&641.3 &2.7&131.6&873.7&2.6&138.4&1408.2 \\
      \multirow{1}{*}{ILU(2)}
       &2.8&87.8&340.2 &2.7&95.7&626.8 &2.7&101.2&859.0 &2.6&105.0&1363.8\\
       \multirow{1}{*}{ILU(4)}
   	&2.8&67.3&351.8&2.7&74.1&650.5 &2.7&78.5&895.2&2.6&81.6&1455.1\\
	\hline
\end{tabular}
\label{tab:ilu-case4}
  }
}
\end{center}
\end{table}

%\begin{table}[htp]
%\begin{center}
%\tabcolsep0.08in
%\extrarowheight=0.7pt
% \caption{The impact of different subdomain solvers for the one-level Schwarz preconditioner for the fixed velocity $u=0.0005 m/s$ case.}
%\begin{tabular}{|c|ccc|ccc|}
%\hline
%End of the simulation (day)&\multicolumn{3}{c|}{t=1} &\multicolumn{3}{c|}{t=2}\\
%\cline{1-4}\cline{5-7}
%Subdomain solver & N. It &L. It &Time & N. It &L. It &Time  \\
%\hline
%
%\hline
%End of the simulation (day)&\multicolumn{3}{c|}{t=3} &\multicolumn{3}{c|}{t=5}\\
%\cline{1-4}\cline{5-7}
%Subdomain solver & N. It &L. It &Time & N. It &L. It &Time  \\
%\hline
%  \multirow{1}{*}{LU}
%       \\
%  \multirow{1}{*}{ILU(0)}
%        \\
%      \multirow{1}{*}{ILU(1)}
%        \\
%      \multirow{1}{*}{ILU(2)}
%       \\
%       \multirow{1}{*}{ILU(4)}
%   	\\
%\hline
%\end{tabular}
%\end{center}
%\end{table}

We then perform test with the 2D and 3D test cases by varying the combination of the overlapping factor $\delta$ and the Schwarz preconditioner types, i.e., the classical-AS \eqref{equ:as-class}, the left-RAS \eqref{equ:as-left}, and the right-RAS \eqref{equ:as-right} preconditioners.  Based on the above observations, we take the choice of the subdomain solver with ILU(1) for the 2D test case and ILU(2) for the 3D test case. The numbers of nonlinear and linear iterations together with the execution time are illustrated in Table \ref{different-overlap}. The results in the table suggest that the more robust combination is the left-RAS preconditioner with the overlapping size $\delta=1$ for the  compromise between the linear iteration and the total computing time. We remark that,  when the overlapping size $\delta=0$,  these preconditioners degenerates into the block-Jacobi preconditioner.

\begin{table}[H]
\begin{center}
 \caption{The impact of different one-level Schwarz preconditioners with several overlapping sizes.
 }
\begin{tabular}{|c|c|ccc|ccc|ccc|}
\hline
Case&& \multicolumn{3}{c|}{Case-1}&\multicolumn{3}{c|}{{Case-2}} &\multicolumn{3}{c|}{Case-4} \\
\cline{1-5}\cline{6-8}\cline{9-11}
Preconditioner&$\delta$& N. It &L. It&Time & {N. It} &{L. It} &{Time}  & N. It &L. It &Time\\
\hline
  \multirow{4}{*}{Classical-AS}
      &0  &3.0   & 80.4 	& 31.5  &{2.8}&{73.9}&{321.6}    &4.2  &53.0  & 211.9  \\
      &1  &3.0   &75.7    	&29.9  &{2.8}&{87.8}&{340.2}  &4.3  &76.1 & 239.6  \\
      &2     &3.0   & 76.9  	&   31.8  &{2.8}&{99.1}&{352.7}   &4.3  &92.6  & 286.7      \\
      &3     &3.0   & 76.8 	&    32.5   &{2.8}&{103.4}& {368.2}  &4.3  &93.8 & 326.6     \\
      \hline
      \multirow{3}{*}{Left-RAS}
   &1	&3.0   &55.8 &28.3  &{2.8}&{43.7}&{306.2} & 4.2  &37.5  &205.0 \\
    &2	&3.0  & 51.9  & 28.8 &{2.8}&{41.6}&{312.8}  &4.3 & 34.0  &218.0    \\
      &3	& 3.0 &49.4  & 28.9  &{2.8}&{40.8}&{322.5}  &4.3  &33.2  &232.3 \\
      \hline
      \multirow{3}{*}{Right-RAS}
      &1	&3.0   & 55.6 &28.2   &{2.8}&{44.5}&{307.5}  & 4.3  &37.1  &207.3 \\
  &2	&3.0  &51.4 &28.0 &{2.8}&{42.0}&{311.4} &4.3  &34.3  &215.8  \\
      &3	&3.0  &48.7 	&28.0  &{2.8}&{41.2}&{326.2}    &4.3  &  33.4  & 223.4 \\
\hline
\end{tabular}
\label{different-overlap}
\end{center}
\end{table}

\subsubsection{Two-level Schwarz preconditioners}
For the reservoir simulation with high accuracy, the supercomputer with a large number of processors  is a must, and therefore  the scalabilities of the algorithm with respect to the number of processors are critically important. As introduced in Section \ref{two-RAS}, in the one-level Schwarz preconditioner,  the average number of linear iterations per Newton iteration grows with the number of processors $N_p$, resulting in the deterioration of the implicit solver. It is clear that some stabilization is needed, which is  achieved by the two-level method with the use of a coarse mesh and the interpolation operators explained in Section \ref{two-RAS}. The performance of the two-level method depends heavy on the two linear solvers defined on the coarse and fine meshes. Here, we restrict ourselves within the framework of Schwarz preconditioned GMRES methods, i.e., we refer to these iterative methods as smoothers on the each levels.

There are several assembly techniques available to construct a hybrid two-level Schwarz preconditioner  by composing the one-level additive Schwarz preconditioner with a coarse-level preconditioner  in a multiplicative or additive manner. Choosing the right type of two-level Schwarz preconditioners is very important for the overall performance of the preconditioner. A large number of numerical experiments is often necessary to identify the right selection. As introduced in Section \ref{two-RAS}, in the study we investigate the pure-coarse \eqref{equ:coarse-level}, the additive \eqref{equ:two-additive}, the left-Kaskade \eqref{equ:left-Kas}, the  right-Kaskade \eqref{equ:right-Kas}, and the  V-cycle type \eqref{equ:v-cycle} two-level Schwarz preconditioners.  For each numerical case, the the overlapping sizes on the fine and coarse levels are fixed to {$\delta_f=1$ and $\delta_c=1$}, respectively. The subdomain solvers on the fine and coarse levels is solved by the LU factorization. In the test, we again run the 2D model with a fixed time step size $0.1$ day {on a $512 \times 512$ mesh} and run Case-4 with a fixed time step size $0.2$ day {on a $128 \times 128 \times 128$ mesh}, and  the simulation is finished at $t= 2$ and $1$ days, respectively. {For the simulation of Case-4,  the flow problem with the fixed velocity model is solved on a $128\times128\times128$ mesh, and the computation is ended at $t=1$ day with $\Delta t=0.1$ day.} The  coarse-to-fine mesh ratio is used to $2$ in each direction. In Table \ref{different-level}, we report the performance of the proposed  two-level  Schwarz preconditioners with respect to different interpolation  operators that includes the first-order scheme \eqref{equ:inter-frist}, the second-order scheme \eqref{equ:inter-second}, and the third-order scheme \eqref{equ:inter-third}.  Below we list the observations made from the results.
\begin{itemize}
\item[(a)] We know that the effectiveness of the Schwarz preconditioner relies on its ability to mimic the spectrum of the linear operator and at the same time is relatively cheap to apply. We see that using the  high order schemes  in the construction of the Schwarz preconditioner provides less linear iteration counts. Moreover,  when compared with the Schwarz method with high order methods,  the  low order approach  is more attractive in the terms of the execution time, owing to its lower bandwidth and a less number of nonzeros in the sparse matrix. The results in the table suggest that the second order scheme is a suitable choice for compromise between the  iterations and the total computing time.
\item[(b)] The best choices for some of the options in the multilevel Schwarz preconditioner  are problem dependent.  For the implicit solution of subsurface flows problems,  from Table \ref{different-level}, we can see that the  additive or  Kaskade type two-level Schwarz preconditioners is exacerbated by a larger number of outer linear iterations, when compared with the V-cycle approach. The pure coarse version of the two-level approach performs considerably worse than the hybrid methods. Hence, our experiments suggest that there is a  benefit to include both pre- and post-swipes of the one-level preconditioning for the simulation of model problems,  especially for the 3D test cases.
\end{itemize}

\begin{table}[hbp]
{\small
\begin{center}
 \caption{The impact of different two-level  Schwarz preconditioners with three interpolation operators. In the table, the symbol ``--" denotes the divergence of the solver caused by the failure of linear iterations.}
\begin{tabular}{|c|l|ccc|ccc|ccc|}
\hline
\multirow{2}{*}{Case}& \multirow{2}{*}{Type} & \multicolumn{3}{c|}{First-order} &\multicolumn{3}{c|}{Second-order}& \multicolumn{3}{c|}{Third-order}\\
\cline{3-5}\cline{6-8}\cline{9-11}
& & N. It &L. It &Time &N. It &L. It &Time & N. It &L. It &Time\\
\hline
\multirow{5}{*}{Case-1}&Coarse
  &--   &--    &-- &-- &-- &-- &-- &-- &--  \\
%\hline
&Additive
  &3.0  &114.6 & 233.9 &3.0 &71.3 &207.3 &3.0 &45.9 &734.7  \\
%\hline
&Left-Kaskade
  &3.0   &71.1 &204.2 &3.0 &34.1 &163.9  &3.0 & 8.05 & 565.2 \\
%\hline
&Right-Kaskade
  &3.0  &75.0  &194.8 &3.0 &37.7 &157.8  &3.0 &  9.03 & 565.5  \\
%\hline
&V-Cycle
  &3.0 &10.0    & 129.6    &3.0    &9.4   &129.1 &3.0  &5.4 &579.3   \\
  \hline
  \multirow{5}{*}{{Case-2}}&{Coarse}
  &{--}   &{--}     &{--}  &{--}  &{--}  &{--}  &{--}   &{--}   &{--}    \\
%\hline
&{Additive}
  &{2.8}&{87.4}&{371.0}&{2.8}&{58.3}&{349.4}&{2.8}&{42.0}&{5186.7}   \\
%\hline
&{Left-Kaskade}
   &{2.8}&{31.3}&{325.7} &{2.8}&{20.6}&{312.4}&{2.8}&{14.7}&{4921.4}   \\
%\hline
&{Right-Kaskade}
   &{2.8}&{34.5}&{326.6}&{2.8}&{21.8}&{315.9}&{2.8}&{15.1}&{4940.3}   \\
%\hline
&{V-Cycle}
  &{2.8}&{14.9}&{314.2}&{2.8}&{12.6} &{308.8}&{2.8}&{5.8}&{4770.9}   \\
\hline
\multirow{5}{*}{Case-4}&Coarse
  &--   &--    &-- &-- &-- &-- &--  &--  &--   \\
%\hline
&Additive
  &4.0    &79.9  &453.2   &4.0  &44.5  &402.1 &4.0  &36.6  &6114.6   \\
%\hline
&Left-Kaskade
  &--   &--  &--  &--  &--  &-- &4.0  &13.4  &5763.8    \\
%\hline
&Right-Kaskade
  &--   &-- &-- &--  &-- &--&4.0  &14.8  &5308.8     \\
%\hline
&V-Cycle
 &4.0   &11.4 &374.2   &4.0  &10.3  &372.8 &4.0  &5.1  &5156.5   \\
\hline
\end{tabular}
\label{different-level}
\end{center}
}
\end{table}

For the two-level preconditioner, the size of the coarse mesh has a strong impact on the efficiency and robustness of  the method. It is clear that using a relatively fine coarse mesh gains a stronger two-level Schwarz preconditioner, and therefore it can  help reduce the total number of linear iterations. On the other hand, finer coarse meshes generates plenty of amount of memory and cache, leading to the increase of the total compute time. An important implementation detail to consider in designing the two-level method is to balance the effects of preconditioning and the computing time of the coarse solve. To understand the impact of different coarse meshes on the convergence of the algorithm, we show the results with different coarse mesh sizes with respect to different interpolation  operators for the 2D and 3D test cases. In the test,  the fine meshes for the 2D and 3D problems are $512 \times 512$ and $128 \times 128 \times 128$, respectively.  As shown in Table \ref{different-interpolation},  we observe from the table that:  (a) for the first-order scheme,  the bad quality of the second coarse mesh does  lead to a large increase of the number of iterations; (b) with the help of higher order schemes, the linear iterations  increase  slowly with growth of  coarse-to-fine mesh ratios. This implies that a higher order  scheme  does a good job on preconditioning the fine mesh problem, since the second coarse mesh even with worse quality is able to keep the whole algorithm efficient  in terms of  the number of iterations.

\begin{table}[hbp]
\begin{center}
 \caption{The impact of coarse-to-fine mesh ratios for the V-Cycle two-level  Schwarz preconditioner. }
\begin{tabular}{|c|c|ccc|ccc|ccc|}
\hline
\multirow{2}{*}{Case}& \multirow{2}{*}{Mesh ratio}& \multicolumn{3}{c|}{First-order} &\multicolumn{3}{c|}{Second-order}& \multicolumn{3}{c|}{Third-order}\\
\cline{3-5}\cline{6-8}\cline{9-11}
&& N. It &L. It&Time & N. It &L. It &Time & N. It &L. It &Time\\
\hline
\multirow{4}{*}{Case-1}
&2  &3.0    & 10.0    	& 129.6    &3.0    &9.4    &129.1 &3.0  &5.4   &579.3  \\
&4  &3.0    & 18.8    	&123.0     &3.0    &16.9    &120.9 &3.0  &13.3   &553.1  \\
&8    &3.0    &26.8     &130.1     &3.0    &19.2    &120.0 &3.0  &14.1   &535.5  \\
&16    &3.0    &36.0    &140.5     &3.0    &19.5    &118.5 &3.0  &14.2  &532.8  \\
\hline
\multirow{4}{*}{{{Case-2}}}
&{2} &{2.8}&{14.9}&{314.2}&{2.8}&{12.6} &{308.8}&{2.8}&{5.8}&{4770.9}   \\
&{4}&{2.8}&{16.2}&{310.7}&{2.8}&{13.1}&{303.1}&{2.8}&{8.9}&{4741.6}\\
&{8}&{2.8}&{18.9} &{308.5}&{2.8}&{13.7}&{299.6}&{2.8}&{9.4}& {4725.3}\\
&{16}&{2.8}&{30.4}&{326.2}&{2.8}&{15.8}&{302.4}&{2.8}&{11.3}&{4709.2}\\
\hline
\multirow{4}{*}{Case-4}
&2  &4.0   &11.4   	&374.2   &4.0   &10.3   &372.8   &4.0  &5.1  &5156.5   \\
&4  &4.0    &17.7   	&372.1   &4.0   & 10.6  &368.4   &4.0  &9.5  &5106.9   \\
&8  &4.0    & 18.8  	&370.2   &4.0   & 11.0  &354.4   &4.0  &10.2  &5091.3   \\
&16 &4.0    &39.8   	&415.7   &4.0   &12.7   &356.3   &4.0  &11.6  &5051.4   \\
\hline
\end{tabular}
\label{different-interpolation}
\end{center}
\end{table}

\subsection{Parallel scalability}
Achieving good parallel scalability is important in parallel computing, especially when solving large-scale problems with many processors.  Hence, in the following we focus on the parallel performance of the proposed fully implicit method with one-level and two-level restricted additive Schwarz preconditioners. Again, the 2D and 3D test problems descried in subsectionn  \ref{sec:Schwraz} are used for the scalability simulations. The numerical tests are carried out on the Tianhe-2 supercomputer. The computing nodes of Tianhe-2 are interconnected via a proprietary high performance network, and there are two 12-core Intel Ivy Bridge Xeon CPUs and 24 GB local memory in each node. In the numerical experiments, we use all 24 CPU cores in each node and assign one subdomain to each core.

\begin{table}[hbp]
\begin{center}
 \caption{Strong scalability with different number of processors $N_p$.}
\begin{tabular}{|c|c|c|ccc|ccc|}
\hline
\multirow{2}{*}{Case}&\multirow{2}{*}{Mesh}& \multirow{2}{*}{$N_p$}& \multicolumn{3}{c|}{One-level}& \multicolumn{3}{c|}{Two-level} \\
\cline{4-6}\cline{7-9}
&&& N. It &L. It&Time & N. It &L. It &Time \\
\hline
%\multirow{12}{*}{Case-1}&\multirow{5}{*}{$4096\times4096$}
%&64   &6.2  &98.1   &545.2  &6.2  &24.5   &409.8 \\
%&&128 &6.2  &136.3   &291.7  &6.2  &25.0   &169.2 \\
%&&256  &6.2  &138.9   &142.7  &6.2  &25.3   &81.9 \\
%&&512  &6.2  &206.5   &88.0 &6.2  &26.1   &38.8 \\
%&&1024  &6.2  &205.1    &43.1  &6.2   &27.0   &20.3 \\
%&&2048  &6.2  &322.8  &33.2  &6.2   &29.8  &13.5 \\
%\cline{2-9}
\multirow{6}{*}{Case-1}&\multirow{6}{*}{$8192\times8192$}
&256   &7.2  &65.4  &564.3 &7.2  & 18.3&460.7  \\
&&512   &7.2  &95.0  &291.5 &7.2  &18.9  &187.6 \\
&&1024   &7.2  &96.3  &136.6 &7.2  &19.6  &90.0 \\
&&2048   &7.2  &142.1  &82.1 &7.2  &20.4  &46.1 \\
&&4096   &7.2  &144.4  &44.9 &7.2  &21.1  &26.3 \\
&&8192  &7.2  &214.2    &34.6  &7.2   &24.0   &18.8
 \\
\hline
%\multirow{10}{*}{Case-4}&\multirow{5}{*}{$128\times128\times128$}
%&64 &3.3 &40.7 &525.7 &3.3 &15.8&504.2 \\
%&&128  &3.3 &47.8   &160.9  &3.3  &16.8    &148.0  \\
%&&256 &3.3 &53.2   &53.6  &3.3  &17.3    &47.8\\
%&&512 &3.3 &60.3    &24.1   &3.3  &17.4    & 19.8  \\
%&&1024 &3.3 &71.4    &12.7  &3.3  &20.1    &10.5   \\
%\cline{2-9}
\multirow{5}{*}{Case-4}&\multirow{5}{*}{$256\times256\times256$}
&512  &7.8 &82.7   &775.7 &7.8    &31.0 &735.1  \\
&&1024  &7.8  &97.7   &250.1    &7.8  &35.2   &231.5  \\
&&2048  &7.8 &109.8    &97.6  &7.8  &37.9  &76.2  \\
&&4096  &7.8 & 122.9   &68.7   &7.8  & 43.2 & 37.6\\
&&8192  &7.8 &147.0    &29.9   &7.8  & 35.0   &21.8  \\
\hline
\end{tabular}
\label{impact-strong-scale}
\end{center}
\end{table}

\begin{figure}[htbp]
\begin{center}
\subfloat[Case-1]{\includegraphics[width=0.45 \textwidth]{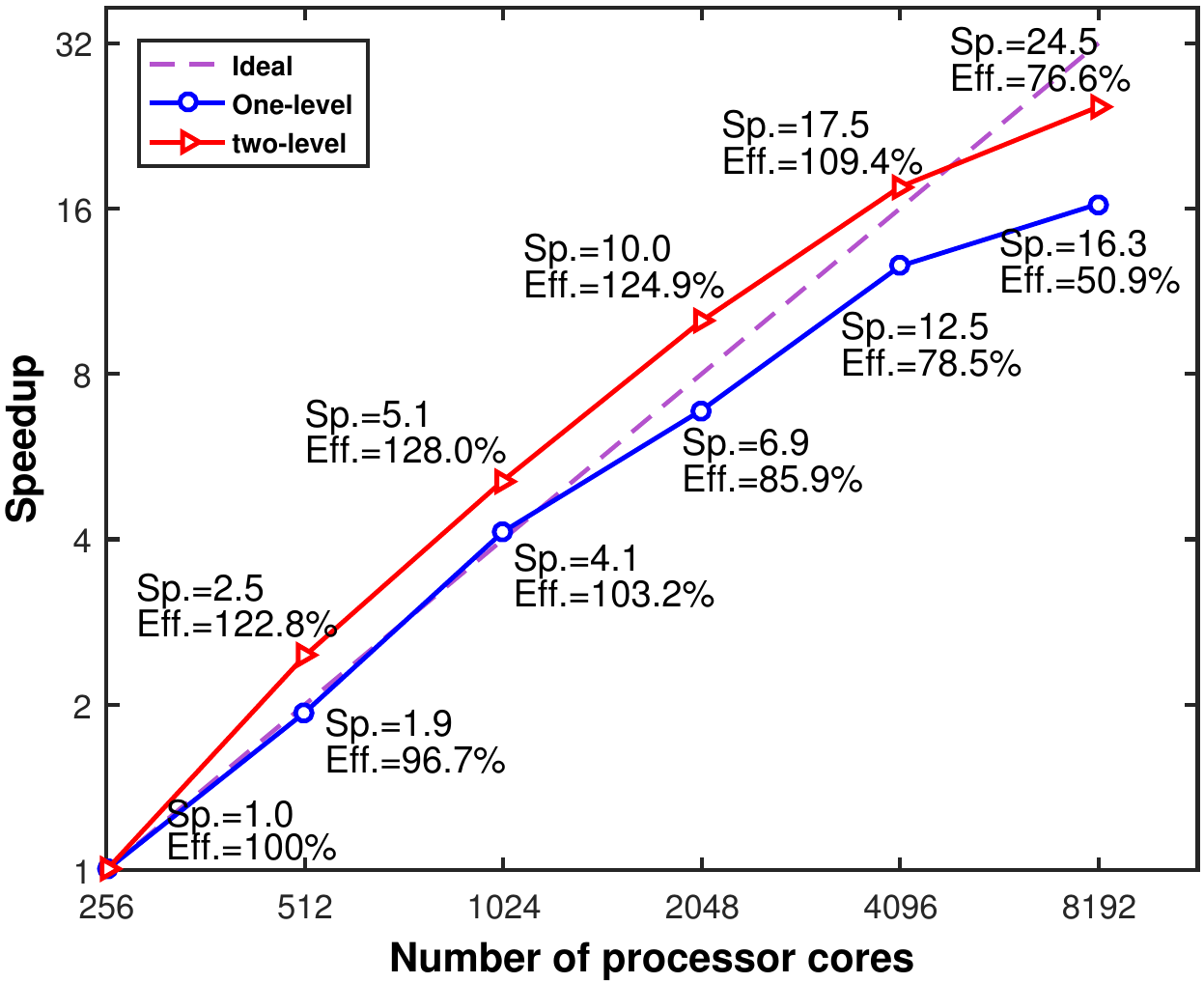}}
\subfloat[Case-4]{\includegraphics[width=0.45 \textwidth]{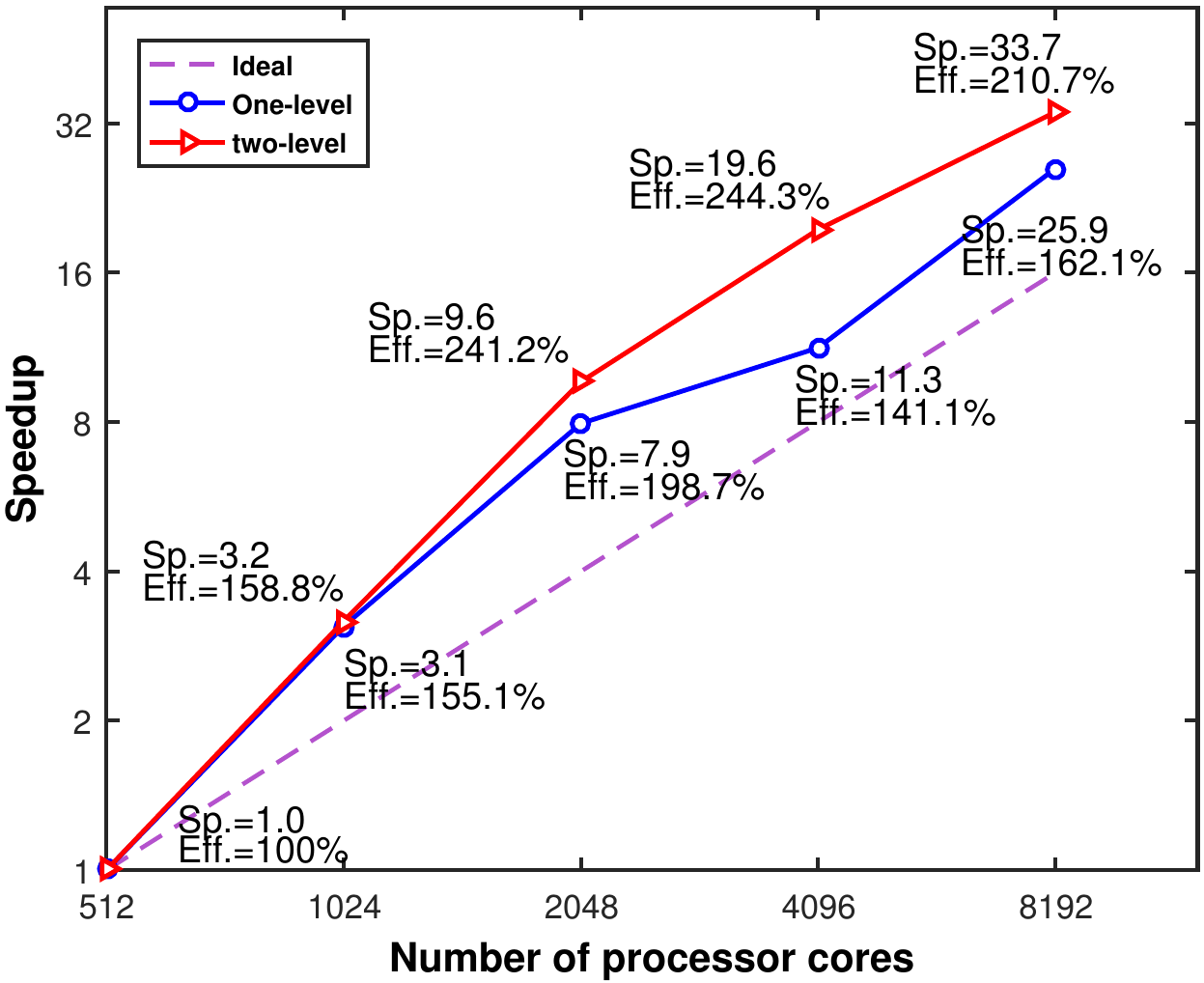}}
\caption{Strong scalability and parallel  efficiency results  with  different number of processors.}
\label{fig:strong-scale}
\end{center}
\end{figure}

The strong scalability  (denoted by $Speedup$) and the parallel efficiency (denoted by $E_{f}$) are respectively defined as follows:
\begin{equation*}
Speedup=\displaystyle\frac{T(N_{min})}{T(N_p)},~E_{f}=\displaystyle\frac{T(N_{min})\times N_{min}}{T(N_p)\times N_p},
\end{equation*}
where $N_{min}$ denotes the smallest processors number of the compassion, $T(N_p)$ denotes computational time {with} $N_p$ processors. In the strong scalability test of the 2D problem, we use a fixed mesh $8192\times 8192$, and  also fixed time step sizes $\triangle t=0.05$ and 0.01 days, in which the largest simulation consists of $8192\times 8192=67,108,864$ degrees of freedom. Moreover, we investigate the scalability of the the proposed fully implicit solver  by using the fixed mesh $256\times 256\times 256$, also a fixed time step size $\triangle t=0.5$ day for Case-4.  The number of processors is changed from $256$ to $8192$.  Table \ref{impact-strong-scale}  shows the number of  nonlinear and linear iterations as well as  the computing time with respect to the number of processors, and  Figure \ref{fig:strong-scale} provides the compute time and $Speedup$ curve for the strong scalability.  The table clearly indicates that the number of Newton iterations remains to be independent of the number of processors, and the number of linear iterations depends on the preconditioner employed in the solver. For the one-level solver, the number of linear iterations suffers as the number of processors increases. While, with the help of the coarse mesh, the number of linear  iterations for the two-level preconditioner  is kept to a low level as the number of processors increases.  {Moreover, to investigate  the strong scalability of the proposed method with respect to the heterogeneity property in the reservoir, we again use the Case-2 benchmark, and compare the number of nonlinear and linear iterations and the total computing time at $t=1,2,3,5$ days, respectively. As shown in Table \ref{tab:strong-scale-case-4}, a good strong scalability is also achieved for the random case test under longer simulations.}

\begin{table}[hbp]
\begin{center}
{
{\footnotesize
 \caption{{Strong scalability with different number of processors $N_p$ for Case-2.}}
\begin{tabular}{|c|c|ccc|ccc|ccc|ccc|}
\hline
\multirow{2}{*}{Level}& \multirow{2}{*}{$N_p$}& \multicolumn{3}{c|}{t=1}& \multicolumn{3}{c|}{t=2} & \multicolumn{3}{c|}{t=3}& \multicolumn{3}{c|}{t=5} \\
\cline{3-14}
&& N. It &L. It&Time & N. It &L. It &Time & N. It &L. It &Time & N. It &L. It &Time\\
\hline
\multirow{5}{*}{One}
&64 &2.8 &44.5  &307.5 &2.7 &52.1 &588.6  &2.7 &58.7 &801.4  &2.6 &63.5&1275.7 \\
&128  &2.8&53.7&125.2  &2.7&64.8&245.2  &2.7&69.3&353.2  &2.6&75.2&583.1\\
&256   &2.8&63.1&57.8  &2.7&71.6&110.9  &2.7&80.1&169.6  &2.6&84.9&252.7\\
&512   &2.8&72.6&28.3  &2.7&81.2&52.6  &2.7&87.7&80.2  &2.6&92.5&134.9\\
&1024   &2.8&84.9&18.9  &2.7&94.5&37.2  &2.7&99.4&58.5  &2.6&104.7&96.3\\
&2048   &2.8&92.5&12.4  &2.7&106.7& 28.5  &2.7&112.8&39.1  &2.6&118.5&62.5\\
\cline{1-14}
\multirow{5}{*}{Two}
&64    &2.8&18.8&283.6  &2.7&20.2&554.7  &2.7&22.6&772.0  &2.6&24.3&1236.8\\
&128   &2.8&21.5&104.5  &2.7&24.5&218.1  &2.7&27.3&312.4  &2.6&29.7&545.9\\
&256    &2.8&25.8&42.6  &2.7&28.9&91.5  &2.7&30.6&137.9  &2.6&33.5&227.5\\
&512   &2.8&29.4&22.8  &2.7&31.7&46.5  &2.7&34.3&74.1  &2.6&36.4&118.4\\
&1024   &2.8&32.3&15.7  &2.7&35.4&32.9  &2.7&36.8&51.8  &2.6&38.2&89.2\\
&2048  &2.8&37.6&11.0 &2.7&39.2&25.3 &2.7&40.2&35.7  &2.6&42.5&58.3\\
\hline
\end{tabular}
\label{tab:strong-scale-case-4}
}
}
\end{center}
\end{table}

The weak scalability is used to examine how the execution time varies with the number of processors when the problem size per processor  is fixed.  In the weak scaling test, we start with a small $2048\times2048$ mesh with the number of processors $N_p=64$  and end up with a large $16384\times16384$ mesh ($268,435,456$ degrees of freedom) using up to $4096$ processor cores for the 2D test case.  Also, we further test our algorithms in terms of the weak scalability starting with a $96\times96\times96$ mesh and $N_p=216$ for Case-4. We refine the mesh and increase the number of processors simultaneously to keep the number of unknowns per processor as a constant.
Table \ref{tab:weak-scale} presents the results of weak scaling tests, which are run with fixed time step sizes $\triangle t=0.01$ and 0.5 days respectively, and are stopped after $5$ implicit time steps, i.e., the simulation are terminated at $t=0.05$ and $2.5$ days.  We observe that, with the increase of the number of processors form $64$ to $4096$ and $216$ to $6859$ for the 2D and 3D test problems, a reasonably good weak scaling performance is obtained, especially for the two-level  restricted Schwarz approach,  which indicates that the proposed solver has a good weak scalability for this range of processor counts. 

\begin{table}[hbpt]
\begin{center}
 \caption{Weak scalability with different number of processors.}
\begin{tabular}{|c|c|c|ccc|ccc|}
\hline
\multirow{2}{*}{Case}&\multirow{2}{*}{Mesh}& \multirow{2}{*}{$N_p$}& \multicolumn{3}{c|}{One-level}& \multicolumn{3}{c|}{Two-level} \\
\cline{4-6}\cline{7-9}
&&& N. It &L. It&Time & N. It &L. It &Time \\
\hline
\multirow{6}{*}{Case-1}
&$2048\times2048$  &64 &3.8  &26.4  &37.4  &3.8  &12.7  &35.2  \\
&$4096\times4096$  &256  &4.8  &46.4  &61.1  &4.8  &18.0  &54.3 \\
&$6144\times6144$   &576     &6.2    &69.5  &95.1   &6.2   &14.4 &67.2  \\
&$8192\times8192$   &1024   &7.2  &96.3  &136.6 &7.2  &19.6  &90.0   \\
&$12288\times12288$ &2304 &9.8 &173.9 &316.8 &9.8 &17.8 &131.8 \\
&$16384\times16384$ &4096 &12.4 &294.7&607.8 &12.4 &17.0  &197.1\\
\hline
\multirow{6}{*}{Case-4}
&$96\times96\times96$  &216  &3.8   & 40.6  &12.5   &3.8  &16.9    &10.7   \\
&$128\times128\times128$  &512  &4.2   &56.3   &15.0    &4.2  &15.9    &12.2   \\
&$192\times192\times192$  &1728  &5.8   &90.3    &27.6   &5.8  &35.9    &22.5   \\
&$224\times224\times224$  &2744  &6.8   &107.7    &35.2   &6.8  &38.2   &28.1   \\
&$256\times256\times256$  &4096   &7.8 & 122.9   &68.7   &7.8  & 43.2 & 37.6\\
&$304\times304\times304$ &6859 &9.4
  &147.4 &85.9  &9.4 &44.7 & 54.6\\
\hline
\end{tabular}
\label{tab:weak-scale}
\end{center}
\end{table}

In summary,  observing from the above tables and the figures, we highlight that, compared with the one-level method,  the two-level restricted additive Schwarz method results in a very sharp reduction in the number of linear iterations and therefore brings about a good reduction in compute time. Hence, the  two-level method  is much more effective and scalable than the one-level approach in terms of the strong and weak scalabilities.

\section{Conclusions}\label{conclusion}
The simulation of subsurface flows with high resolution solutions is of paramount importance in reservoir  simulation.
In this work, we have presented a parallel fully implicit framework based on multilevel restricted Schwarz preconditioners for subsurface flow simulations with Peng-Robinson equation of state. The proposed framework can get rid of  the restriction of the time step size, and  is flexible and allows us to construct different type of multilevel preconditioners, based on plenitudinous choices for additively or multiplicatively strategies, interpolation and restriction operators.    After experimenting with many different overlapping Schwarz type preconditioners,  we found that the class of V-cycle-type restricted Schwarz methods based on the second  order scheme is extremely beneficial for the problems under investigation. Numerical experiments also showed that the proposed algorithms and simulators are robust  and scalable for the large-scale solution of some benchmarks as well as realistic problems with highly heterogeneous permeabilities in petroleum reservoirs.

\section*{Acknowledgment}
\red{The authors would like to express their appreciations to the anonymous reviewers for the invaluable comments that have greatly improved the quality of the manuscript.}

\end{document}